\documentclass[10pt,twocolumn,twoside]{IEEEtran}
\ifCLASSINFOpdf
 
\else
  
\fi

\usepackage{amsmath,amssymb,amsthm}
\usepackage{color}
\usepackage{balance}
\usepackage[noadjust]{cite}
\usepackage{graphicx}
\usepackage{epsfig}
\usepackage{breqn}
\usepackage{cite}
\usepackage[caption=false,font=footnotesize]{subfig}
\usepackage{array}
\usepackage{scalerel}

\newcommand{\norm}[1]{\left\lVert#1\right\rVert}

\newcommand{\mycolor} {black}

\newtheorem{theorem}{Theorem}[section]

\newtheorem{lemma}[theorem]{Lemma}
\newtheorem{prop}[theorem]{Proposition}
\newtheorem{cor}[theorem]{Corollary}

\allowdisplaybreaks
\newcommand\given[1][]{\:#1\vert\:}
\renewcommand\thesection{\arabic{section}}
\newcommand{\bP}[1]{{\mathbb{P}}\left[{#1}\right]}
\newcommand{\bE}[1]{{\mathbb{E}}\left[{#1}\right]}

\newcommand{\1}[1]{{\bf 1}\left[#1\right]}       

\newcommand{\fsquare}{\vrule height6pt width7pt depth1pt}   
\newcommand{\myproof}{{\hfill \\ \bf Proof. \ }}           
\newcommand{\myendpf}{\hfill\fsquare \\[0.1in]}             

\begin{document}
\bstctlcite{IEEEexample:BSTcontrol}

\title{Connectivity of Wireless Sensor Networks Secured by Heterogeneous Key Predistribution Under an On/Off Channel Model}

\author{Rashad~Eletreby,~\IEEEmembership{Student Member,~IEEE,}
        and~Osman~Ya\u{g}an,~\IEEEmembership{Senior Member,~IEEE,}

\thanks{R. Eletreby and O. Ya\u{g}an are with the Department
of Electrical and Computer Engineering and CyLab, Carnegie Mellon University, Pittsburgh,
PA, 15213 USA. E-mail: reletreby@cmu.edu, oyagan@ece.cmu.edu.}
\thanks{This work has been supported in part by
National Science Foundation through grant CCF \#1617934 and by a generous gift from Persistent Systems, Inc.
}
}

\maketitle

\begin{abstract}

We investigate the connectivity of a wireless sensor network secured by the heterogeneous key predistribution scheme under an independent on/off channel model. The heterogeneous scheme induces an inhomogeneous random key graph, denoted by $\mathbb{K}(n;\pmb{\mu},\pmb{K},P)$ and the on/off channel model induces an Erd\H{o}s-R\'enyi graph, denoted by $\mathbb{H}(n,\alpha)$. Hence, the overall random graph modeling the WSN is obtained by the intersection of $\mathbb{K}(n;\pmb{\mu},\pmb{K},P)$ and $\mathbb{H}(n,\alpha)$. We present conditions on how to scale the parameters of the intersecting graph with respect to the network size $n$ such that the graph i) has no isolated nodes and ii) is connected, both with high probability as the number of nodes gets large. Our results are supported by a simulation study demonstrating that i) despite their asymptotic nature, our results can in fact be useful in designing {\em finite}-node wireless sensor networks so that they achieve secure connectivity with high probability; and ii) despite the simplicity of the on/off communication model, the probability of connectivity in the resulting wireless sensor network approximates very well the case where the disk model is used. 
\end{abstract}


\IEEEpeerreviewmaketitle

\vspace{-4mm}
\section{Introduction}

\subsection{Wireless Sensor Networks and Security}
\IEEEPARstart{W}{ireless} sensor networks (WSNs) emerged as an enabling platform for a broad range of application areas owing to their low-cost, low-power, small size, and adaptability to the physical environment \cite{Akyildiz_2002}. These unique features triggered the proliferation and adoption of WSNs in several domains including military, health, and environment, but also gave rise to unique security challenges that can not be tackled using classical security mechanisms \cite{security_survey}. In particular, {\em asymmetric} cryptosystems provide a scalable solution for securing large scale WSNs; however, they are generally slow and {\color{\mycolor} lead to} excessive energy and memory consumption. On the other hand, {\em symmetric} cryptosystems were shown to be superior in terms of speed and energy efficiency, but they demand novel and efficient mechanisms for key-establishment among sensor nodes \cite{Gligor_2002,Haowen_2003}. In principle, an efficient key-establishment mechanism should result in a {\em {\color{\mycolor} securely} connected} topology, i.e., a network where there exists a secure communication path (possibly multihop) between every pair of nodes allowing the exchange of data and control messages, while conforming to the typical limitations of WSNs. Also, it shall not assume knowledge of post-deployment configuration, {\color{\mycolor} since} in most cases WSNs are deployed randomly in large {\color{\mycolor} numbers}.

In their seminal work, Eschenauer and Gligor proposed {\color{\mycolor} a {\em random key predistribution} protocol} as a practical and efficient method for key-establishment in large scale WSNs \cite{Gligor_2002}. {\color{\mycolor} Their} scheme, hereafter referred to as the EG scheme, operates as follows: before deployment, each node is given a {\em random} set of $K$ cryptographic keys, selected {\em uniformly} (without replacement) from a large key pool of size $P$. After deployment, two nodes can communicate {\em securely} over an existing channel {\em if} they share at least one key. The EG scheme is currently regarded as one of the most feasible solutions for key-establishment among sensor nodes, e.g., see \cite[Chapter~13]{Raghavendra_2004}, \cite{camtepe_2005}, and references therein, and has led the way to several other variants, including the $q$-composite scheme \cite{Haowen_2003}, the random pairwise scheme \cite{Haowen_2003}, and many others.

The EG scheme inherently assumes that all nodes are {\em homogeneous} in terms of their roles and capabilities, hence they are assigned the same number $K$ of keys. However, emerging WSN applications are complex and are envisioned to require the coexistence of different classes of nodes with different roles and capabilities \cite{Yarvis_2005}. For instance, a particular class of nodes may act as cluster heads that are used to connect several clusters of nodes together. These cluster heads need to communicate with a large number of nodes in their vicinity and they are also expected to be more powerful than regular nodes. Thus, one can reasonably argue that more keys should be given to the cluster heads to ensure high levels of connectivity and security.

To cope with the expected heterogeneity in WSN topologies, Ya\u{g}an proposed a new variation of the EG scheme, referred to as the {\em heterogeneous} random key predistribution scheme \cite{Yagan/Inhomogeneous}. The heterogeneous scheme considers the case when the network includes sensors with varying levels of resources, features, security, or connectivity requirements. The scheme is described as follows. Given $r$ classes, each sensor is independently classified as a class-$i$ node with probability $\mu_i>0$ 
for each $i=1,\ldots, r$. Then, sensors in class-$i$ are each assigned $K_i$ keys selected uniformly at random  from a key pool of size $P$. Similar to the EG scheme, nodes that share at least one common key (regardless of their class) can communicate securely over an available channel after deployment.

Given the randomness involved in the EG scheme and the heterogeneous scheme, there is a positive probability that a pair of nodes may have no common key, thus can not establish a secure communication link in between. Moreover, two nodes that share a key may not have a wireless channel in between (possibly because of the limited transmission radius). Hence, it is natural to ask whether the resulting network would be securely {\em connected} or not. Specifically, two nodes are securely connected if they share a key {\em and} have a communication channel in between. A network is said to be connected if there is a path between every pair of vertices. In essence, one {\color{\mycolor} needs} to know if it is possible to {\em control} the parameters of the scheme (possibly as functions of the network size $n$), such that the resulting network is connected {\color{\mycolor} with high probability}. 

In \cite{Yagan/Inhomogeneous}, Ya\u{g}an considered a WSN secured by the heterogeneous scheme under {\em full-visibility} assumption, i.e., all pairs of sensors have a communication channel in between, hence the only condition for two nodes to be connected is to share a key. Therein, {\color{\mycolor} they} established scaling conditions on the parameters of the heterogeneous scheme as functions of the network size $n$ such that the resulting network is connected with high probability as the number of nodes gets large. In particular, {\color{\mycolor} they} considered a random graph model naturally induced by the heterogeneous scheme and established scaling conditions on the model parameters such that the resulting graph is connected with high probability as the number of nodes gets large. Specifically, with $\pmb{K}=\{K_1,K_2,\ldots,K_r \}$, $\pmb{\mu}=\{\mu_1,\mu_2,\ldots,\mu_r \}$, and $n$ denoting the network size, we let $\mathbb{K}(n; \pmb{\mu},\pmb{K},P)$ denote the random graph induced by the heterogeneous key predistribution scheme, where any pair of vertices are {\em adjacent} as long as they share a key. This model was referred to as the {\em inhomogeneous} random key graph in \cite{Yagan/Inhomogeneous}, where zero-one laws for absence of isolated nodes and connectivity were established. The inhomogeneous random key graph models the {\em shared-key} connectivity of the WSN under the heterogeneous scheme.

Our paper is motivated by the fact that the full-visibility assumption is not likely to hold in real-world implementations of WSNs. In particular, the randomness of the wireless channel as well as limited transmission ranges would severely limit the availability of wireless channels between nodes, rendering two nodes disconnected even when they share a key. In fact, as wireless connectivity comes into play, an essential question {\color{\mycolor} arises:} {\em Under a given model for wireless connectivity, is it possible to control the parameters of the heterogeneous scheme to ensure that the resulting network is connected?}

\vspace{-3mm}
\subsection{Modeling Wireless Connectivity}

Our paper aims to {\color{\mycolor} answer} this question, hence bridging the disconnect between the model developed in \cite{Yagan/Inhomogeneous} and real world implementations of WSNs where wireless channels are scarce and the full-visibility assumption does not hold. In particular, we model the {\em wireless connectivity} of the WSN, say using a (possibly random) graph $\mathbb{I}(n;\cdot)$, whose edges represent pairs of sensors who have a wireless communication channel available in between. The overall model of the WSN will then be an intersection of $\mathbb{K}(n; \pmb{\mu},\pmb{K},P)$ and $\mathbb{I}(n;\cdot)$ since a pair of sensors can establish a {\em secure communication link} if they share a key {\em and} have a wireless channel available. Let $\mathbb{G}$ be the intersecting graph, i.e., $\mathbb{G} := \mathbb{K}(n; \pmb{\mu},\pmb{K},P) \cap \mathbb{I}(n;\cdot)$. At a high level, our objective is to establish scaling conditions on the parameters of $\mathbb{G}$ such that the resulting graph is connected with high probability as the number of nodes gets large. 

In practice, limited transmission {\color{\mycolor} range of sensors significantly impacts} the wireless connectivity of a WSN, hence the {\em disk model} \cite{Gupta99} can be seen as a good candidate model for wireless connectivity among sensor nodes. The disk model is described as follows. Assuming that nodes are distributed over a bounded region $\mathcal{D}$ of a euclidean plane, nodes $v_i$ and $v_j$ located at $\pmb{x}_i$ and $\pmb{x}_j$, respectively, are able to communicate if $\norm{\pmb{x}_i-\pmb{x}_j} < \rho$, where $\rho$ denotes the transmission radius. A special case of the disk model when node locations are independently and uniformly distributed over the region $\mathcal{D}$, gives rise to the random geometric graph \cite{PenroseBook}, hereafter denoted $\mathbb{I}(n;\rho)$. Now, let $\mathbb{G}(n; \pmb{\mu},\pmb{K},P,\rho)$ be a random graph obtained by intersecting the inhomogeneous random key graph $\mathbb{K}(n; \pmb{\mu},\pmb{K},P,\cdot)$ with a random geometric graph $\mathbb{I}(n;\rho)$. Clearly, $\mathbb{G}(n; \pmb{\mu},\pmb{K},P,\rho)$ represents {\color{\mycolor} a reasonably accurate} model for a WSN secured by the heterogeneity scheme, where two nodes are connected if they i) share a key, and ii) are within transmission radius. 

Unfortunately, analyzing the connectivity of  $\mathbb{G}(n; \pmb{\mu},\pmb{K},P,\rho) $ is likely to be very challenging, and may very well be impossible. In fact, the Gupta-Kumar conjecture \cite{Gupta99} on the connectivity of $\mathbb{H}(n;\alpha) \cap \mathbb{I}(n;\rho)$ where $\mathbb{H}(n;\alpha)$ represents an Erd\H{o}s-R\'enyi (ER) graph, took {\color{\mycolor} many} years (and several attempts) to be resolved {\color{\mycolor} eventually} by Penrose \cite{penrose2016connectivity}; see \cite{Yagan/EG_intersecting_ER} for a detailed discussion on the difficulties involved in analyzing {\em intersection} of different types of graphs. {\color{\mycolor} The model $\mathbb{K}(n; \pmb{\mu},\pmb{K},P)$ considered here} is much more complicated than an ER graph {\color{\mycolor} due to edge correlations \cite{Yagan/Inhomogeneous}, leading to the following important question:} {\em Is there any communication model that provides a good approximation of the classical disk model, but also allows a comprehensive analysis of the resulting intersecting graph?}

This question was {\color{\mycolor} answered in the affirmative in \cite{Yagan/EG_intersecting_ER,YaganPairwise}}, where it was shown that an independent on/off channel model provides a good approximation of the disk model {\color{\mycolor} for understanding the critical scalings of connectivity} in settings similar to ones we consider here. In the independent on/off channel model, the wireless channel between any given pair of nodes is either on (with probability $\alpha$) or off (with probability $1-\alpha$) independently from all other channels. The model induces an ER graph $\mathbb{H}(n;\alpha)$, where an edge exists (respectively does not exist) between two vertices with probability $\alpha$ (respectively $1-\alpha$) independently from all other edges.

{\color{\mycolor} With these in mind}, we model the wireless connectivity of the WSN by an ER graph $\mathbb{H}(n;\alpha)$ and study the connectivity of the intersecting graph $\mathbb{G}(n; \pmb{\mu},\pmb{K},P,\alpha):=\mathbb{K}(n; \pmb{\mu},\pmb{K},P) \cap \mathbb{H}(n;\alpha)$. This approach allows us to i) establish rigorous results concerning the secure connectivity of a WSN albeit using a simplified wireless communication model, and ii) demonstrate via simulations that these results still apply under the more realistic disk model. In Section~\ref{sec:numerical}, we provide simulation results indicating that the connectivity of $\mathbb{K}(n; \pmb{\mu},\pmb{K},P) \cap \mathbb{I}(n;\rho)$ behaves very similar to that of $\mathbb{K}(n; \pmb{\mu},\pmb{K},P) \cap \mathbb{H}(n;\alpha)$, as we {\em match} $\alpha$ and $\rho$ leading to the same probability of wireless channel availability; i.e., $\alpha = \pi \rho^2$.
\newline

\vspace{-3mm}
\subsection{Contributions}

We investigate the connectivity of a WSN secured by the heterogeneous key predistribution scheme under an independent on/off channel model. The heterogeneous scheme induces an inhomogeneous random key graph, denoted by $\mathbb{K}(n;\pmb{\mu},\pmb{K},P)$ and the on/off channel model induces an ER graph, denoted by $\mathbb{H}(n,\alpha)$. Hence, the overall random graph modeling the WSN is obtained by the intersection of $\mathbb{K}(n;\pmb{\mu},\pmb{K},P)$ and $\mathbb{H}(n,\alpha)$. We denote this intersection by $\mathbb{G}(n;\pmb{\mu},\pmb{K},P,\alpha)$, i.e., $\mathbb{G}(n;\pmb{\mu},\pmb{K},P,\alpha):=\mathbb{K}(n;\pmb{\mu},\pmb{K},P) \cap \mathbb{H}(n,\alpha)$. We present conditions on how to scale the parameters of $\mathbb{G}(n;\pmb{\mu},\pmb{K},P,\alpha)$ with respect to the network size $n$ such that i) it has no isolated nodes and ii) it is connected, both with high probability as the number of nodes gets large. {\color{\mycolor} The results are given in the form of zero-one laws with critical scalings precisely established}. This maps to {\em dimensioning} the parameters of the heterogeneous scheme with respect to the network size $n$ and the channel parameter $\alpha$ such that the resulting network is securely connected.

Our results are supported by a simulation study (see Section \ref{sec:numerical}) demonstrating that i) despite their asymptotic nature, our results can in fact be useful in designing {\em finite}-node WSNs so that they achieve secure connectivity with high probability; and ii) despite the simplicity of the on/off communication model, the probability of connectivity in the resulting WSN approximates very well the case where the disk model is used. In addition, our results are shown to complement and generalize several previous work in the literature (see Section~\ref{sec:results} for details). 

All limiting statements, including asymptotic equivalences are considered with the number of sensor nodes $n$ going to infinity. The indicator function of an event $E$ is denoted by $\pmb{1}[E]$. 
We say that an event holds with high probability (whp) if it holds with probability $1$ as $n \rightarrow \infty$. 
In comparing the asymptotic behavior of the sequences $\{a_n\},\{b_n\}$, we use the standard Landau notation, e.g.,
$a_n = o(b_n)$,  $a_n=\omega(b_n)$, $a_n = O(b_n)$, $a_n = \Omega(b_n)$, and
$a_n = \Theta(b_n)$. We also use $a_n \sim b_n$ to denote the asymptotic equivalence $\lim_{n \to \infty} {a_n}/{b_n}=1$.

\section{The Model}
The heterogeneous random key predistribution scheme introduced
in \cite{Yagan/Inhomogeneous} works as follows. Consider a network of $n$ sensors labeled as $v_1, v_2, \ldots,v_n$. Each sensor node is classified into one of the $r$ classes (e.g., priority levels) according to a probability distribution $\pmb{\mu}=\{\mu_1,\mu_2,\ldots,\mu_r\}$ with $\mu_i >0$ for $i=1,\ldots,r$ and $\sum_{i=1}^r \mu_i=1$. Then, a class-$i$ node is assigned $K_i$ cryptographic keys selected uniformly at random from a key pool of size $P$. It follows that the key ring $\Sigma_x$ of node $v_x$ is a random variable (rv) with 
\begin{equation} \nonumber
\mathbb{P}[\Sigma_x=S \mid t_x=i]= \binom P{K_i}^{-1}, \quad S \in \mathcal{P}_{K_i},
\nonumber
\end{equation}
where $t_x$ denotes the class of $v_x$ and $\mathcal{P}_{K_i}$ is the collection of
all subsets of $\{1,\ldots, P\}$ with size $K_i$.
The classical key predistribution scheme of Eschenauer and Gligor \cite{Gligor_2002}
constitutes a special case of this model with $r=1$, i.e., when all sensors belong to the same class and receive the same number of keys; see also \cite{yagan2012zero}.

Let $\pmb{K}=\{K_1,K_2,\ldots,K_r\}$ and assume without loss of generality that $K_1 \leq K_2 \leq \ldots \leq K_r$. Consider a random graph $\mathbb{K}$ induced on the vertex set $\mathcal{V}=\{v_1,\ldots,v_n\}$ such that a pair of nodes $v_x$ and $v_y$ are adjacent, denoted by $v_x \sim_{\mathbb{K}} v_y$, if they have at least one cryptographic key in common, i.e.,
\begin{equation}
v_x \sim_{\mathbb{K}} v_y \quad \text{if} \quad \Sigma_x \cap \Sigma_y \neq \emptyset.
\label{adjacency_condition}
\end{equation}

The adjacency condition (\ref{adjacency_condition}) defines the inhomogeneous random key graph denoted by $\mathbb{K}(n;\pmb{\mu},\pmb{K},P)$  \cite{Yagan/Inhomogeneous}. This model is also known in the literature  as
the {\em general random intersection graph}; e.g., see \cite{Zhao_2014,Rybarczyk,Godehardt_2003}. 
 The probability $p_{ij}$ that a class-$i$ node and a class-$j$ node are adjacent is given by
\begin{equation}
p_{ij} = \mathbb{P}[v_x \sim_{\mathbb{K}} v_y \mid t_x=i, t_y=j] =
1-\frac{\binom {P-K_i}{K_j}}{\binom {P}{K_j}}
\label{eq:osy_edge_prob_type_ij}
\end{equation}
as long as $K_i + K_j \leq P$; otherwise if $K_i +K_j > P$, we  have $p_{ij}=1$.
Let $\lambda_i$ denote the \textit{mean} probability that a class-$i$ node is connected to another node in $\mathbb{K}(n;\pmb{\mu},\pmb{K},P)$. We have
\begin{align}
\lambda_i & =\mathbb{P}[v_x \sim_{\mathbb{K}} v_y \mid t_x=i ]
=\sum_{j=1}^r p_{ij} \mu_j.
 \label{eq:osy_mean_edge_prob_in_RKG}
\end{align}
We also find it useful to define the {\em mean} key ring size by $K_{\textrm{avg}}$; i.e.,
\begin{equation}
K_{\textrm{avg}} = \sum_{j=1}^{r} K_j \mu_j.
\label{eq:mean_key_size}
\end{equation}

We model the wireless connectivity of the WSN by means of an independent on/off channel model. In particular, the channel between any given pair of nodes is either on with probability $\alpha$ or off with probability $1-\alpha$. More precisely, let $\{B_{ij}(\alpha), 1 \leq i < j \leq n\}$ denote i.i.d Bernoulli rvs, each with success probability $\alpha$. The communication channel between two distinct nodes $v_x$ and $v_y$ is on (respectively, off) if $B_{xy}(\alpha)=1$ (respectively, if $B_{xy}(\alpha)=0$). The on/off channel model induces a standard ER graph $\mathbb{H}(n;\alpha)$  \cite{Bollobas}, defined  on the vertices $\mathcal{V}=\{v_1,\ldots,v_n\}$ such that $v_x$ and $v_y$ are adjacent, denoted by $v_x \sim_{\mathbb{H}} v_y$, if $B_{xy}(\alpha)=1$.

We model the overall topology of a WSN by the intersection of an inhomogeneous random key graph $\mathbb{K}(n;\pmb{\mu},\pmb{K},P)$ with an ER graph $\mathbb{H}(n;\alpha)$. Namely, nodes  $v_x$ and $v_y$ are adjacent in $\mathbb{K} (n;\pmb{\mu},\pmb{K},P) \cap \mathbb{H}(n;\alpha)$, if and only if they are adjacent in both $\mathbb{K}$ \textit{and} $\mathbb{H}$. Hence, the edges in the intersection graph
$\mathbb{K} (n;\pmb{\mu},\pmb{K},P) \cap \mathbb{H}(n;\alpha)$ represent pairs of sensors that
can securely communicate since they have i) a communication link available in between, and ii) a shared cryptographic key.
Therefore, studying the connectivity properties of $\mathbb{K} (n;\pmb{\mu},\pmb{K},P) \cap \mathbb{H}(n;\alpha)$ amounts to studying the secure connectivity of heterogeneous WSNs under the on/off channel model.

Hereafter, we denote the intersection graph $\mathbb{K} (n;\pmb{\mu},\pmb{K},P) \cap \mathbb{H}(n;\alpha)$ by $\mathbb{G}(n;\pmb{\mu},\pmb{K},P,\alpha)$. To simplify the notation, we let $\pmb{\theta}=(\pmb{K},P)$, and $\pmb{\Theta}=(\pmb{\theta},\alpha)$. The probability of edge existence between a class-$i$ node $v_x$ and a class-$j$ node $v_y$ in $\mathbb{G}(n;\pmb{\Theta})$ is given by
\begin{align} 
&\mathbb{P}[v_x \sim_{\mathbb{G}} v_y \given[\Big] t_x=i,t_y=j] \nonumber \\
&=\mathbb{P}[v_x \sim_{\mathbb{K}} v_y \cap v_x \sim_{\mathbb{H}} v_y \given[\big] t_x=i,t_y=j] \nonumber \\
&=\alpha p_{ij} \nonumber
\end{align}
by independence. Similar to (\ref{eq:osy_mean_edge_prob_in_RKG}), the mean edge probability for a class-$i$ node in $\mathbb{G}(n;\pmb{\mu},\pmb{\Theta})$ as $\Lambda_i$ is given by
\begin{align} 
\Lambda_i = \sum_{j=1}^r \mu_j \alpha p_{ij} = \alpha \lambda_i, \quad i=1,\ldots, r.
\label{eq:osy_mean_edge_prob_in_system}
\end{align}

Throughout, we assume that the number of classes $r$ is fixed and does not scale with $n$, and so are the probabilities $\mu_1, \ldots,\mu_r$. All of the remaining parameters are assumed to be scaled with $n$.

\vspace{-3mm}
\section{Main Results and Discussion}
We refer to a mapping $\pmb{\Theta}=K_1,\ldots,K_r,P,\alpha:
 \mathbb{N}_0 \rightarrow \mathbb{N}_0^{r+1} \times (0,1)$ as a \textit{scaling} 
if
\begin{equation}
1 \leq K_{1,n} \leq K_{2,n} \leq \ldots \leq K_{r,n} \leq P_n/2
\label{scaling_condition_K}
\end{equation}
for all $n=2,3,\ldots$. We note that under (\ref{scaling_condition_K}), the edge probability $p_{ij}$ is given by
(\ref{eq:osy_edge_prob_type_ij}).

\subsection{Results}
\label{sec:results}
We first present a zero-one law for the absence of isolated nodes in 
$\mathbb{ G}(n;\pmb{\mu},\pmb{\Theta}_n)$.
\begin{theorem}
\label{theorem:isolated_nodes}
{\sl
Consider a probability distribution $\pmb{\mu}=\{\mu_1,\ldots,\mu_r\}$ with $\mu_i >0$ for $i=1,\ldots,r$ and a scaling 
$\pmb{\Theta}: \mathbb{N}_0 \rightarrow \mathbb{N}_0^{r+1} \times (0,1)$ such that
\begin{equation}
\Lambda_1(n)=\alpha_n \lambda_1(n) \sim c \frac{\log n}{n}
\label{scaling_condition_KG}
\end{equation}
for some $c>0$. We have
\begin{equation}
\lim_{n\to\infty} \mathbb{P} \left[\begin{split} \mathbb{G}(n;\pmb{\mu},\pmb{\Theta}_n) \text{ has~} \\ \text{ no isolated nodes} \end{split} \right ]=
\begin{cases} 
      \hfill 0    \hfill & \text{ if $c<1$} \\
      \hfill 1 \hfill & \text{ if $c>1$} \\
\end{cases}
\label{eq:node_isolation}
\end{equation}
}
\end{theorem}
The scaling condition (\ref{scaling_condition_KG}) will often be used in the form
\begin{equation} \label{scaling_condition_KG_v2}
\Lambda_1(n)=c_n \frac{\log n}{n}, \ n=2,3,\ldots
\end{equation}
with $\lim_{n\to\infty} c_n=c>0$. 

Next, we present an analogous result for connectivity.
\begin{theorem}
\label{theorem:connectivitiy}
{\sl
Consider a probability distribution $\pmb{\mu}=\{\mu_1,\ldots,\mu_r\}$ with $\mu_i >0$ for $i=1,\ldots,r$ and a scaling $\pmb{\Theta}: \mathbb{N}_0 \rightarrow \mathbb{N}_0^{r+1} \times (0,1)$ such that
(\ref{scaling_condition_KG})
holds for some $c>0$. Then, we have
\begin{equation} 
\lim_{n\to\infty} \mathbb{P}[\mathbb{ G}(n;\pmb{\mu},\pmb{\Theta}_n) \text{ is connected}]=
\begin{cases} 
      \hfill 0    \hfill & \text{ if $c<1$} \\
      \hfill 1 \hfill & \text{ if $c>1$} \\
\end{cases}
\label{eq:connectivity}
\end{equation}
under the additional conditions that
\begin{equation}
P_n \geq \sigma n, \quad n=1, 2, \ldots
\label{eq:conn_Pn2}
\end{equation}
for some $\sigma>0$
and
\begin{equation}
p_{11}(n) =\omega\left(\frac{1}{n \alpha_n}\right).
\label{eq:conn_K1}
\end{equation}
}
\end{theorem}

The resemblance of the results presented in Theorem~\ref{theorem:isolated_nodes} and Theorem~\ref{theorem:connectivitiy} indicates that  absence of isolated nodes and connectivity are asymptotically equivalent properties for $\mathbb{ G}(n;\pmb{\mu},\pmb{\Theta}_n)$. Similar observations were made
for other well-known random graph models as well; e.g., inhomogeneous random key graphs \cite{Yagan/Inhomogeneous}, Erd\H{o}s-R\'enyi graphs \cite{Bollobas}, and (homogeneous) random key graphs \cite{yagan2012zero}.

Conditions (\ref{eq:conn_Pn2}) and (\ref{eq:conn_K1}) are enforced mainly for technical reasons and they are only needed in the proof of the one-law of Theorem~\ref{theorem:connectivitiy}. In particular, condition (\ref{eq:conn_Pn2}) is essential for real-world WSN implementations in order to ensure the {\em resilience} of the network against node capture attacks; e.g., see \cite{Gligor_2002,DiPietroTissec}. For instance, assume that an adversary captures a number of sensors, compromising all the keys that belong to the captured nodes. If $P_n = o(n)$, then it would be possible for the adversary to compromise $\Omega(P_n)$ keys by capturing only $o(n)$ sensors (whose type does not matter in this case). In this case, the WSN would fail to exhibit the {\em unassailability} property \cite{MeiPanconesiRadhakrishnan2008,YM_ToN} and would be deemed as vulnerable against adversarial attacks. 

Also, condition (\ref{eq:conn_K1}) is enforced mainly for technical reasons for the proof of the one-law to work. The need of such a lower bound arises from the fact that our scaling condition (\ref{scaling_condition_KG}) merely scales the minimum {\em mean} edge probability, not the minimum (or each) edge probability, as $\log n / n$. For instance, the current scaling condition (\ref{scaling_condition_KG}) gives us an easy upper bound on the minimum edge probability in the network, but does not specify any non-trivial lower bound on that probability. More specifically, it is easy to see that $\alpha_n p_{11}(n) = O \left( \Lambda_1 \right) = O \left( \log n / n \right)$, but it is not clear if the sequence $\alpha_n p_{11}(n)$ has a non-trivial lower bound. In fact, authors in \cite{devroye2014connectivity} investigated the connectivity of an inhomogeneous Erd\H{o}s-R\'enyi (ER) graph, while setting the probability of an edge connecting two nodes of classes $i$ and $j$ to $\kappa \left(i,j \right) \log n / n$, where $\kappa \left(i,j \right)$ returns a positive real number for each pair $(i,j)$; i.e., each individual edge was scaled as $\log n / n$. 

In summary, condition (\ref{eq:conn_Pn2}) is needed to ensure the resilience of the network against node capture attacks, while condition (\ref{eq:conn_K1}) is needed to provide a non-trivial lower bound on the minimum edge probability of the network. {\color{\mycolor} To provide a concrete example, one can set $P_n = n \log n$ and have $K_{1,n} =  (\log n)^{1/2+\varepsilon}$  with any $\varepsilon > 0$ to satisfy (\ref{eq:conn_K1}) for any $\alpha_n \geq 1/(\log n)^{\varepsilon}$ (see Lemma \ref{lemma:conn_asym_eq}).
In this case, setting $K_{\textrm{avg},n}=\log n^{3/2}$ ensures that the resulting network is connected whp (see Corollary \ref{cor:osy_Kav}).}

Theorem~\ref{theorem:isolated_nodes} (resp. Theorem~\ref{theorem:connectivitiy}) states that $\mathbb{ G}(n;\pmb{\mu},\pmb{\Theta}_n)$ has no isolated node (resp. is connected) whp if the mean degree of class-$1$ nodes (that receive the smallest number $K_{1,n}$ of keys) is scaled as $(1+\epsilon) \log n$ for some $\epsilon > 0$. On the other hand, if this minimal mean degree scales as $(1-\epsilon) \log n$ for some $\epsilon > 0$, then whp $\mathbb{G}(n;\pmb{\mu},\pmb{\Theta}_n)$ has an isolated node, and hence not connected. 
These results indicate that the minimum key ring size in the network has a 
significant impact on the {\color{\mycolor} connectivity of 
$\mathbb{G}(n;\pmb{\mu},\pmb{\Theta}_n)$.}

The importance of the minimum key ring size on connectivity can be seen more explicitly under a mild condition on the scaling, as shown in the next corollary.

\vspace{-1mm}
\begin{cor}
\label{cor:osy_Kav}
{\sl Consider a probability distribution $\pmb{\mu}=\{\mu_1,\ldots,\mu_r\}$ with $\mu_i >0$ for $i=1,\ldots,r$ and a scaling $\pmb{\Theta}: \mathbb{N}_0 \rightarrow \mathbb{N}_0^{r+1} \times (0,1)$ such that $\lambda_1(n)=o(1)$ and
\begin{equation}
\alpha_n \frac{K_{1,n} K_{\textrm{avg},n}}{P_n} \sim c \frac{\log n}{n}
\label{eq:osy_new_scaling_Kav}
\end{equation}
for some $c > 0$, where $K_{\textrm{avg},n}$ is as defined at (\ref{eq:mean_key_size}).
Then we have the zero-one law (\ref{eq:node_isolation}) for absence of isolated nodes. If, in addition, the conditions (\ref{eq:conn_Pn2}) and (\ref{eq:conn_K1}) are satisfied, then we also have the zero-one law (\ref{eq:connectivity}) for connectivity.
}
\end{cor}

\vspace{-4mm}
\myproof
In view of  (\ref{eq:osy_mean_edge_prob_in_RKG}), we see that $\lambda_1(n)=o(1)$ implies $p_{1j}(n)=o(1)$ for $j=1,\ldots,r$. From Lemma~\ref{lemma:conn_asym_eq}, this then leads to $p_{1j}(n) \sim \frac{K_{1,n}K_{j,n}}{P_n}$, whence
\[
\lambda_1(n) = \sum_{j=1}^r \mu_j p_{1j}(n) \sim \frac{K_{1,n} \sum_{j=1}^r \mu_jK_{j,n}} {P_n}  =\frac{K_{1,n} K_{\textrm{avg},n}}{P_n} 
\]
Thus, the scaling conditions (\ref{scaling_condition_KG}) and (\ref{eq:osy_new_scaling_Kav}) are equivalent under $\lambda_1(n)=o(1)$
and Corollary \ref{cor:osy_Kav} follows from Theorem~\ref{theorem:isolated_nodes} and Theorem~\ref{theorem:connectivitiy}.
\myendpf
We see from Corollary~\ref{cor:osy_Kav} that for a fixed mean number $K_{\textrm{avg},n}$ of keys per sensor, network connectivity is directly affected by the minimum key ring size $K_{1,n}$. For example, reducing $K_{1,n}$ by half means that the smallest $\alpha_n$ for which the network becomes connected whp is increased by two-fold (see Figure \ref{fig:2} for a numerical example demonstrating this phenomenon).

\vspace{-3mm}
\subsection{Comparison with related work}
\label{sec:comparison}

Our main results extend the work in \cite{Yagan/Inhomogeneous} and \cite{Jun_2016}, where authors established zero-one laws for the connectivity of a WSN secured by the heterogeneous key predistribution scheme under the {\em full-visibility} assumption. Although a crucial first step in the study of heterogeneous key predistribution schemes, the assumption that all pairs of sensors have a communication channel in between is not likely to hold in most practical settings. In this regard, our work extends the results in \cite{Yagan/Inhomogeneous} and \cite{Jun_2016} to more practical WSN scenarios where the wireless connectivity of the network is taken into account. By setting $\alpha_n=1$ for each $n=1,2,\ldots$ (i.e., by assuming that all links are available), our results reduce to those given in \cite{Yagan/Inhomogeneous}. 

Authors in \cite{Yagan/EG_intersecting_ER} (respectively, \cite{zhaoconnectivity}) investigated the connectivity (respectively, $k$-connectivity) of WSNs secured by the classical EG scheme under an independent on/off channel model. However, when the network consists of sensors with varying level of resources (e.g., computational, memory, power), and with varying level of security and connectivity requirements, it may no longer be sensible to assign the same number of keys to all sensors. Our work addresses this issue by generalizing  \cite{Yagan/EG_intersecting_ER} to the cases where nodes can be assigned different number of keys.  When $r=1$, i.e., when all nodes belong to the same class and receive the same number of keys, our result recovers the main result in  \cite{Yagan/EG_intersecting_ER}.

\section{Numerical Results}
\label{sec:numerical}
We now present numerical results to support Theorems \ref{theorem:isolated_nodes} and Theorem~\ref{theorem:connectivitiy} in the finite node regime. Furthermore, we show by simulations that the on/off channel model serves as a good approximation of the disk model. In our simulations, we fix the number of nodes at $n = 500$ and the size of the key pool at $P = 10^4$. 

The first step in comparing the on/off channel model to the disk model is to propose a {\em matching} between ER graph $\mathbb{H}(n;\alpha)$ and the random geometric graph $\mathbb{I}(n;\rho)$ in a way that leads to the same probability of link availability. In particular, consider $500$ nodes distributed uniformly and independently over a folded unit square $[0,1]^2$ with toroidal (continuous) boundary conditions. Since there are no border effects, we get
\begin{equation} \nonumber
\mathbb{P}\left[\norm{\pmb{x}_i-\pmb{x}_j} < \rho \right] = \pi \rho^2, \quad \quad i \neq j, \quad i,j=1,\ldots,n
\end{equation}
whenever $\rho < 0.5$. Thus, in order to match the two communication models we set $\alpha=\pi \rho^2$. Recall that $\mathbb{G}(n;\pmb{\mu},\pmb{\theta},\alpha)=\mathbb{K}(n;\pmb{\mu},\pmb{\theta}) \cap \mathbb{H}(n;\alpha)$, and let $\mathbb{\widetilde{G}}(n;\pmb{\mu},\pmb{\theta},\rho)=\mathbb{K}(n;\pmb{\mu},\pmb{\theta}) \cap \mathbb{I}(n;\rho)$. Next, we present several simulation results comparing the (empirical) probabilities that $\mathbb{G}$ and $\mathbb{\widetilde{G}}$ are connected, respectively. 

We start by considering the channel parameter $\alpha = \pi \rho^2 = 0.2$, $\alpha = \pi \rho^2 =0.4$, $\alpha = \pi \rho^2 =0.6$, and $\alpha = \pi \rho^2 = 0.8$, while varying the parameter $K_1$ (i.e., the smallest key ring size) from $5$ to $35$. The number of classes is fixed at $2$ with $\pmb{\mu}=\{0.5,0.5\}$ and we set $K_2=K_1+5$. For each parameter pair $(\pmb{K}, \alpha)$ (respectively, $(\pmb{K}, \pi \rho^2)$), we generate $800$ independent samples of the graphs $\mathbb{G}$ (respectively, $\mathbb{\widetilde{G}}$) and count the number of times (out of a possible $800$) that the obtained graphs i) have no isolated nodes and ii) are connected. Dividing the counts by $800$, we obtain the (empirical) probabilities for the events of interest.  We observed that $\mathbb{G}$ is connected whenever it has no isolated nodes yielding the same empirical probability for both events. This is in parallel with the asymptotic equivalence of the two properties as implied by Theorems \ref{theorem:isolated_nodes} and \ref{theorem:connectivitiy}.

\begin{figure}[t]
\centerline{\includegraphics[scale=0.45]{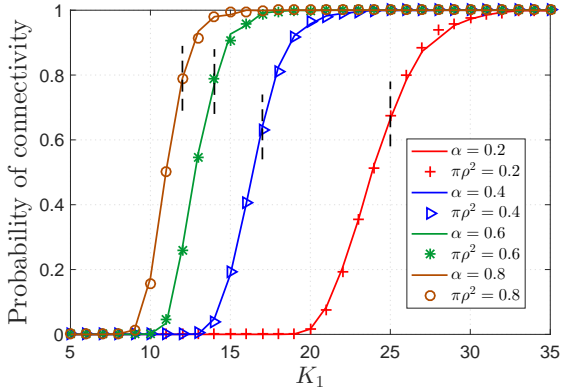}}
\caption{Empirical probability that $\mathbb{G}$ and $\mathbb{\widetilde{G}}$ are connected as a function of $\pmb{K}$ for $\alpha = \pi \rho^2=0.2$, $\alpha = \pi \rho^2=0.4$, $\alpha = \pi \rho^2= 0.6$, and $\alpha = \pi \rho^2=0.8$ with $n = 500$ and $P = 10^4$; in each case, the empirical probability value is obtained by averaging over $800$ experiments. Vertical dashed lines stand for the critical threshold of connectivity asserted by Theorem~\ref{theorem:connectivitiy}.}
\label{fig:1}
\vspace{-4mm}
\end{figure}

In Figure~\ref{fig:1}, we show the empirical probabilities of the connectivity of $\mathbb{G}$ (represented by lines) and $\mathbb{\widetilde{G}}$ (represented by symbols). We observe that the empirical probabilities are almost identical, supporting the claim that the on/off channel model serves as a good approximation of the disk model under the given matching condition. Furthermore, we show the critical threshold of connectivity {\em predicted} by Theorem~\ref{theorem:connectivitiy} by a vertical dashed line for each curve. More specifically, for a given $\alpha$, the vertical dashed lines stand for the minimum integer value of $K_1$ that satisfies
\begin{equation}
\hspace{-.6mm}\lambda_1(n)\hspace{-1mm}=\hspace{-1mm}\sum_{j=1}^2 \mu_j \hspace{-1mm}\left( \hspace{-.2mm} 1- \frac{\binom{P-K_j}{K_1}}{\binom{P}{K_1}}  \hspace{-.2mm}\right) \hspace{-.1mm}>  \hspace{-.1mm}\frac{1}{\alpha} \frac{\log n}{n}
\label{eq:numerical_critical}
\end{equation}
According to Theorem~\ref{theorem:connectivitiy}, at this critical value of $K_1$ the network would be connected with probability $1$ as the number of nodes {\em tends to infinity}. We see from Figure~\ref{fig:1} that even in the finite-node regime ($n=500$), the critical value of $K_1$ results in a connected network with high probability.

\begin{figure}[t]
\centerline{\includegraphics[scale=0.45]{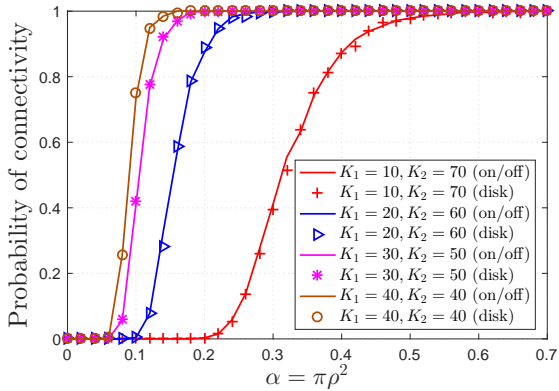}}
\caption{Empirical probability that $\mathbb{G}$ and $\mathbb{\widetilde{G}}$ are connected as a function of $\alpha$ and $\pi \rho^2$ for four choices of $\pmb{K}=(K_1,K_2)$, each with the same mean.}
\label{fig:2}
\vspace{-5mm}
\end{figure}

Figure~\ref{fig:2} is generated in a similar manner with Figure~\ref{fig:1}, this time with an eye towards understanding the impact of the minimum key ring size $K_1$ on network connectivity. We fix the number of classes at $2$ with $\pmb{\mu}=\{0.5,0.5\}$ and consider 
four different key ring sizes $\pmb{K}$
each with mean $40$; we consider
$\pmb{K} = \{10,70\}$, $\pmb{K} = \{20,50\}$, $\pmb{K} = \{30,50\}$, and $\pmb{K} = \{40,40\}$.
We compare the probability of connectivity in the resulting networks as $\alpha$ (respectively, $\pi \rho^2$) varies from zero to one. Although the average number of keys per sensor is kept constant in all four cases, network connectivity improves dramatically as the minimum key ring size $K_1$ increases; e.g., with $\alpha=\pi \rho^2=0.2$, the probability of connectivity is one when $K_1=K_2=40$ while it drops to zero if we set $K_1=10$ and $K_2=70$ so that the mean key ring size is still 40. This confirms the observations made via Corollary \ref{cor:osy_Kav}.

\vspace{-3mm}
\section{Other application areas: The spread of epidemics and information in real-world social networks }

The last decade has witnessed a tremendous advance in our understanding of how information \cite{Yagan_IP2013,Bakshy_2012}, influence \cite{Yagan_Influence2012,Kempe_2003}, {\color{\mycolor} and} diseases \cite{Eubank_2004,Newman_2002} propagate across the globe. A large variety of mathematical models as well as a multitude of data sets paved the way for precise predictions and control of the behavior of such spreading processes on complex networks. In particular, several generative models were proposed to create networks which resemble the structure of real-world complex networks, allowing for large-scale simulations and precise predictions of how a spreading process would behave in real-life. Three structural properties in particular, the {\em power-law degree distribution}, {\em small-world}, and {\em clustering} were shown to be prevalent in real-world social networks \cite{Barabasi_1999,Watts_1998,Amaral_2000}. 

The {\em homogeneous} random key graph (where all nodes receive the same number $K$ of objects) was shown to generate networks that are clustered and small-world \cite{Yagan2017_SmallWorld}. Indeed, the inhomogeneous counterpart $\mathbb{K}(n;\pmb{\mu},\pmb{K},P)$ intrinsically exhibits these two properties as well. In addition to that, one can tune the parameters of $\mathbb{K}(n;\pmb{\mu},\pmb{K},P)$ to generate networks with a power-law degree distribution similar to that observed in real-world social networks \cite{deijfen_kets_2009}. Collectively, the inhomogeneous random key graph $\mathbb{K}(n;\pmb{\mu},\pmb{K},P)$ generates networks {\color{\mycolor} that are small-world and have {\em tunable} degree distribution and clustering, hence it can be considered as a useful model for real-world social networks. In fact, the inhomogeneous random key graph is a natural model for {\em common-interest} social networks.} A common interest relationship between two friends manifests from their selection of common interests {\color{\mycolor} or hobbies from a large pool} \cite{zhaoconnectivity}. Clearly, this can be modeled by an inhomogeneous random key graph, where each individual has a set of interests (possibly of different sizes) sampled from a large pool of interests and two individuals are connected if they happen to share an interest.

{\color{\mycolor} In addition, the intersection model $\mathbb{G}(n;\pmb{\mu},\pmb{K},P,\alpha)$ considered here can be useful in studying the propagation of epidemics or information on complex networks.} A simple model for the spread of epidemics (or information) on complex networks is the so called Susceptible-Infected-Recovered (SIR) model. Therein, a disease is transmitted to a susceptible individual upon contact with an infected individual. Later on, infected individuals recover from the disease and gain immunity from it. The outbreak size is precisely the number of recovered individuals at the steady state. This model results in reasonable predictions for the cases where recovery grants lasting resistance. In \cite{Newman_2002}, it was shown that under some conditions, the dynamics of the SIR model on a given network maps to a {\em bond-percolation} problem with the average {\em transmissibility} of the disease as the percolation parameter. Namely, with $\alpha$ being the average transmissibility; If we are to occupy each edge in the graph with probability $\alpha$, the final outbreak size would be the size of the cluster of vertices that can be reached from the initial infected vertex by traversing the occupied edges only \cite{Newman_2002}. Typically, {\color{\mycolor} one is} interested in deriving the threshold value of $\alpha$ for which a {\color{\mycolor}{\em giant}} connected component emerges, indicating that the disease has reached a positive fraction of the population. 

Intersecting the inhomogeneous random key graph $\mathbb{K}(n;\pmb{\mu},\pmb{K},P)$ with an ER graph $\mathbb{H}(n;\alpha)$ is essentially equivalent to {\em occupying} each edge of $\mathbb{K}(n;\pmb{\mu},\pmb{K},P)$ independently with probability $\alpha$. Hence, the scaling condition for which the one-law of Theorem~\ref{theorem:connectivitiy} holds gives us a threshold value of $\alpha$ for which a strain of a disease or a piece of information would infect the entire population. In particular, let $\hat{\alpha}_n := \log n / \left( n \lambda_1(n) \right)$; if the average transmissibility of a disease $\alpha$ satisfies $\alpha_n>\hat{\alpha}_n$, a single giant component containing all of the vertices emerge (because in this case the network is connected by virtue of Theorem~\ref{theorem:connectivitiy}), allowing the disease to infect each single vertex. Therefore, our results on the connectivity of $\mathbb{G}(n;\pmb{\mu},\pmb{K},P,\alpha)$ {\color{\mycolor} provide} a threshold on the average transmissibility a disease should have (possibly through evolution) in order to persist in a given population modeled by the inhomogeneous random key graph.

\vspace{-3mm}
\section{Proof of Theorem~\ref{theorem:isolated_nodes}}
\label{sec:proof_isolated}

\subsection{Preliminaries}
Few technical results are collected here for convenience. A full list of preliminaries is given in Appendix~\ref{app:prelim}. The first result follows easily from the scaling condition (\ref{scaling_condition_K}). 

\begin{prop} [{\cite[Proposition~4.1]{Yagan/Inhomogeneous}}]
For any scaling $K_1,K_2,\ldots,K_r,P:\mathbb{N}_0 \rightarrow \mathbb{N}_0^{r+1}$, we have (in view of (\ref{scaling_condition_K}))
\begin{equation}
\lambda_1(n) \leq \lambda_2(n) \leq \ldots \leq \lambda_r(n)
\label{eq:isolated_ordering_of_lambda}
\end{equation}
for each $n=2,3,\ldots$.
\end{prop}

Another useful bound that will be used throughout is 
\begin{equation}
(1 \pm x) \leq e^{\pm x}, \quad x \in (0,1)
\label{eq:isolated:exp_bound}
\end{equation}

Finally, we find it useful to write
\begin{equation}
\log (1-x)=-x-\Psi(x)
\label{eq:isolated_log_decomp}
\end{equation}
where  $\Psi(x)=\int_{0}^{x} \frac{t}{1-t} \ \text{dt}$.
From L'H\^{o}pital's Rule, we have
\begin{equation}
\lim_{x\to 0}  \frac{\Psi(x)}{x^2}=\frac{-x-\log (1-x)}{x^2}=\frac{1}{2}.
\label{eq:isolated_hopital}
\end{equation}

\vspace{-3mm}
\subsection{Establishing the one-law}

The proof of Theorem~\ref{theorem:isolated_nodes} relies on the method of first and second moments applied to the number of isolated nodes in $\mathbb{G}(n;\pmb{\mu},\pmb{\Theta}_n)$. Let $I_n(\pmb{\mu},\pmb{\Theta}_n)$ denote the total number of isolated nodes in $\mathbb{G}(n;\pmb{\mu},\pmb{\Theta}_n)$, namely,
\begin{equation}
I_n(\pmb{\mu},\pmb{\Theta}_n)=\sum_{\ell=1}^n \pmb{1}[v_\ell \text{ is isolated in }\mathbb{G}(n;\pmb{\mu},\pmb{\Theta}_n)]
\label{eq:isolated_In}
\end{equation}

The method of first moment \cite[Eqn. (3.10), p.
55]{JansonLuczakRucinski} gives
\begin{equation} \nonumber
1-\mathbb{E}[I_n(\pmb{\mu},\pmb{\Theta}_n)]\leq \mathbb{P}[I_n(\pmb{\mu},\pmb{\Theta}_n)=0] 
\end{equation}

It is clear that in order to establish the one-law, namely that $ \lim_{n \to \infty} \mathbb{P}\left[  I_n(\pmb{\mu},\pmb{\Theta_n})=0\right]=1$, we need to show that
\begin{equation}
\lim_{n \to \infty} \mathbb{E}[I_n(\pmb{\mu},\pmb{\Theta}_n)]=0.
\label{eq:node_isol_one_law_to_show}
\end{equation}

Recalling (\ref{eq:isolated_In}), we have
\begin{align}
&\mathbb{E}\left[I_n(\pmb{\mu},\pmb{\Theta}_n)\right] \nonumber \\
&=n \sum_{i=1}^r \mu_i \mathbb{P}\left[v_1 \text{ is isolated in }\mathbb{G}(n;\pmb{\mu},\pmb{\Theta}_n) \given[\big] t_1=i\right]\nonumber \\
&=n \sum_{i=1}^r \mu_i \mathbb{P}\left[\cap_{j=2}^n [v_j \nsim v_1] \mid v_1 \text{ is class i}\right]\nonumber \\
&=n \sum_{i=1}^r \mu_i \left(\mathbb{P}\left[v_2 \nsim v_1 \mid v_1 \text{ is class i}\right]\right)^{n-1} \label{eq:isolated_independence}
\end{align}
where (\ref{eq:isolated_independence}) follows by the independence of the rvs $\{v_j \nsim v_1\}_{j=1}^n$ given $\Sigma_1$. By conditioning on the class of $v_2$, we find
\begin{align}
\mathbb{P}[v_2 \nsim v_1 \given[\big] t_1=i]
&=\sum_{j=1}^r \mu_j \mathbb{P}[v_2 \nsim v_1 \given[\big] t_1=i,t_2=j]\nonumber \\
&=\sum_{j=1}^r \mu_j (1-\alpha p_{ij}) =1-\Lambda_i \label{eq:isolated_OneLaw_first_step}
\end{align}
Using (\ref{eq:isolated_OneLaw_first_step}) in (\ref{eq:isolated_independence}), and recalling (\ref{eq:isolated_ordering_of_lambda}) and (\ref{eq:isolated:exp_bound}), we obtain
\begin{align*}
\mathbb{E}[I_n(\pmb{\mu},\pmb{\Theta}_n)] &= n \sum_{i=1}^r \mu_i \left(1-\Lambda_i(n)\right)^{n-1}\nonumber \\
&\leq n \left(1-\Lambda_1(n)\right)^{n-1} \leq  e^{\log n \left(1-c_n  \frac{n-1}{n}\right)}.
\end{align*}
Taking the limit as $n$ goes to infinity, we immediately get
(\ref{eq:node_isol_one_law_to_show})
since $\lim_{n \to \infty} (1-c_n  \frac{n-1}{n})=1-c < 0$ under the enforced assumptions (with $c>1$) and the one-law is established.
\myendpf

\vspace{-10mm}
\subsection{Establishing the zero-law}
Our approach in establishing the zero-law relies on the method of second moment applied to a variable that counts the number of nodes that are class-$1$ and isolated. Clearly if we can show that whp there exists at least one class-$1$ node that is isolated under the enforced assumptions (with $c<1$) then the zero-law would immediately follow.

Let $Y_n(\pmb{\mu},\pmb{\Theta}_n)$ denote the number of nodes that are class-$1$ and isolated in $\mathbb{G}(n;\pmb{\mu},\pmb{\Theta}_n)$, and let
\begin{equation} \nonumber
x_{n,i}(\pmb{\mu},\pmb{\Theta}_n)=\pmb{1}[t_i=1 \cap v_i \text{ is isolated in }\mathbb{G}(n;\pmb{\mu},\pmb{\Theta}_n)],
\end{equation}
then we have $Y_n(\pmb{\mu},\pmb{\Theta}_n)=\sum_{i=1}^n x_{n,i}(\pmb{\mu},\pmb{\Theta}_n)$. By applying the method of second moments
\cite[Remark 3.1, p. 55]{JansonLuczakRucinski}  on $Y_n(\pmb{\mu},\pmb{\Theta}_n)$, we get
\begin{equation}
\mathbb{P}[Y_n(\pmb{\mu},\pmb{\Theta}_n)=0] \leq 1-\frac{\mathbb{E}[Y_n(\pmb{\mu},\pmb{\Theta}_n)]^2}{\mathbb{E}[Y_n(\pmb{\mu},\pmb{\Theta}_n)^2]}  
\label{eq:isolated_ZeroLaw_bound}
\end{equation}
where
\begin{equation}
\mathbb{E}[Y_n(\pmb{\mu},\pmb{\Theta}_n)]=n \mathbb{E}[x_{n,1}(\pmb{\mu},\pmb{\Theta}_n)]
\label{eq:isolated_ZeroLaw_first_part}
\end{equation}
and
\begin{align}
\mathbb{E}[Y_n(\pmb{\mu},\pmb{\Theta}_n)^2]=&n \mathbb{E}[x_{n,1}(\pmb{\mu},\pmb{\Theta}_n)] \label{eq:isolated_ZeroLaw_second_part} \\
&+n(n-1)\mathbb{E}[x_{n,1}(\pmb{\mu},\pmb{\Theta}_n) x_{n,2}(\pmb{\mu},\pmb{\Theta}_n)] \nonumber
\end{align}
by exchangeability and the binary nature of the rvs $\{x_{n,i}(\pmb{\mu},\pmb{\Theta}_n) \}_{i=1}^n$.
Using (\ref{eq:isolated_ZeroLaw_first_part}) and  (\ref{eq:isolated_ZeroLaw_second_part}), we get
\begin{equation} \nonumber
\begin{split}
\frac{\mathbb{E}[Y_n(\pmb{\mu},\pmb{\Theta}_n)^2]}{\mathbb{E}[Y_n(\pmb{\mu},\pmb{\Theta}_n)]^2}  =& \frac{1}{n \mathbb{E}[x_{n,1}(\pmb{\mu},\pmb{\Theta}_n)]} \\
&+ {\frac{n-1}{n} \frac{\mathbb{E}[x_{n,1}(\pmb{\mu},\pmb{\Theta}_n) x_{n,2}(\pmb{\mu},\pmb{\Theta}_n)]}{\mathbb{E}[x_{n,1}(\pmb{\mu},\pmb{\Theta}_n)]^2}}
\end{split}
\end{equation}

In order to establish the zero-law, we need to show that
\begin{equation}
\lim_{n \to \infty} n \mathbb{E}[x_{n,1}(\pmb{\mu},\pmb{\Theta}_n)]= \infty,
\label{eq:isolated_ZeroLaw_first_condition}
\end{equation}
and
\begin{equation}
\limsup_{n \to \infty} \left(\frac{\mathbb{E}[x_{n,1}(\pmb{\mu},\pmb{\Theta}_n) x_{n,2}(\pmb{\mu},\pmb{\Theta}_n)]}{\mathbb{E}[x_{n,1}(\pmb{\mu},\pmb{\Theta}_n)]^2}\right) \leq 1.
\label{eq:isolated_ZeroLaw_second_condition}
\end{equation}

{\prop
\label{prop:prop6.1}
Consider a scaling $K_1,\ldots,K_r,P:\mathbb{N}_0 \rightarrow \mathbb{N}_0^{r+1}$ and a scaling $\alpha:\mathbb{N}_0 \rightarrow (0,1)$ such that (\ref{scaling_condition_KG}) holds with $\lim_{n \to \infty} c_n=c>0$. Then, we have
\begin{equation*}
\lim_{n \to \infty} n \mathbb{E}[x_{n,1}(\pmb{\mu},\pmb{\Theta}_n)]= \infty, \quad \text{if } c<1
\end{equation*}
}

\vspace{-5mm}
\begin{proof}
We have
\begin{align}
&n \mathbb{E}\left[x_{n,1}(\pmb{\mu},\pmb{\Theta}_n)\right] \nonumber \\
&=n  \mathbb{P}\left[v_1 \text{ is isolated in }\mathbb{G}(n;\pmb{\mu},\pmb{\Theta}_n) \cap t_1=1 \right]\nonumber \\
&=n \mu_1 \mathbb{P}\left[\cap_{j=2}^n [v_j \nsim v_1] \given[\big] t_1=1\right]\nonumber \\
&=n \mu_1 \mathbb{P}\left[v_2 \nsim v_1\given[\big] t_1=1\right]^{n-1}\nonumber \\
&=n \mu_1 \left(\sum_{j=1}^r \mu_j \mathbb{P}\left[v_2 \nsim v_1 \given[\big] t_1=1,t_2=j\right]\right)^{n-1}\nonumber \\
&=n \mu_1 \left(\sum_{j=1}^r \mu_j (1-\alpha_n p_{1j})\right)^{n-1} \label{eq:int_isol_prob_osy} \\
&=n \mu_1 \left(1-\Lambda_1(n)\right)^{n-1} = \mu_1 e^{\beta_n} 
\label{eq:isolated_ZeroLaw_simp1}
\end{align}
where $\beta_n=\log n+(n-1)\log (1-\Lambda_1(n))$. Recalling (\ref{eq:isolated_log_decomp}), we get
\begin{align}
\beta_n&=\log n-(n-1)\left(\Lambda_1(n)+\Psi(\Lambda_1(n))\right)\nonumber \\
&=\log n-(n-1)\left(c_n \frac{\log n}{n}+\Psi \left(c_n \frac{\log n}{n}\right)\right)\nonumber \\
&=\log n \left(1-c_n \frac{n-1}{n}\right) \nonumber \\
& \quad -(n-1) \left(c_n \frac{\log n}{n} \right)^2 \frac{\Psi \left(c_n \frac{\log n}{n}\right)}{\left(c_n \frac{\log n}{n}\right)^2}
 \label{eq:isolated_ZeroLaw_simp2}
\end{align}

Recalling (\ref{eq:isolated_hopital}), we have
\begin{equation}
\lim_{n \to \infty} \frac{\Psi \left(c_n \frac{\log n}{n}\right)}{\left(c_n \frac{\log n}{n}\right)^2} = \frac{1}{2}
\label{eq:isolated_ZeroLaw_simp3}
\end{equation}
since $c_n \frac{\log n}{n}=o(1)$. Thus, $\beta_n=\log n \left(1-c_n \frac{n-1}{n}\right)-o(1)$.
Using (\ref{eq:isolated_ZeroLaw_simp1}), (\ref{eq:isolated_ZeroLaw_simp2}), (\ref{eq:isolated_ZeroLaw_simp3}), and letting $n$ go to infinity, we get
\begin{equation*}
\lim_{n \to \infty} n \mathbb{E}[x_{n,1}(\pmb{\mu},\pmb{\Theta}_n)]= \infty
\end{equation*}
whenever $\lim_{n \to \infty} c_n=c < 1$.
\end{proof}

{\prop
Consider a scaling $K_1,\ldots,K_r,P:\mathbb{N}_0 \rightarrow \mathbb{N}_0^{r+1}$ and a scaling $\alpha:\mathbb{N}_0 \rightarrow (0,1)$ such that (\ref{scaling_condition_KG}) holds with $\lim_{n \to \infty} c_n=c>0$. Then, we have (\ref{eq:isolated_ZeroLaw_second_condition}) if $c<1$.
\label{prop:new_osy}
}

The proof of Proposition~\ref{prop:new_osy} is given in Appendix~\ref{app:proof_of_prop_6.2}. Collectively, Proposition~\ref{prop:prop6.1} and Proposition~\ref{prop:new_osy} establish (\ref{eq:isolated_ZeroLaw_first_condition}) and (\ref{eq:isolated_ZeroLaw_second_condition}) respectively, which in turn establish the zero-law of Theorem~\ref{theorem:isolated_nodes}.

\vspace{-4mm}
\section{Proof of Theorem~\ref{theorem:connectivitiy}}
Let $C_n(\pmb{\mu},\pmb{\Theta}_n)$ denote the event that the graph $\mathbb{G}(n,\pmb{\mu},\pmb{\Theta}_n)$ is connected, and with a slight abuse of notation, let $I_n(\pmb{\mu},\pmb{\Theta}_n)$ denote the event that the graph $\mathbb{G}(n,\pmb{\mu},\pmb{\Theta}_n)$ has no isolated nodes. It is clear that if a random graph is connected then it does not have any isolated node, hence
\begin{equation*}
C_n(\pmb{\mu},\pmb{\Theta}_n) \subseteq I_n(\pmb{\mu},\pmb{\Theta}_n)
\end{equation*}
and we get
\begin{equation}
\mathbb{P}[C_n(\pmb{\mu},\pmb{\Theta}_n)] \leq \mathbb{P}[I_n(\pmb{\mu},\pmb{\Theta}_n)]
\label{eq:conn_ZeroLaw}
\end{equation}

\vspace{-3mm}
and
\begin{align} \label{eq:conn_OneLaw}
\mathbb{P}[C_n(\pmb{\mu},\pmb{\Theta}_n)^c] & = \mathbb{P}[I_n(\pmb{\mu},\pmb{\Theta}_n)^c]+\mathbb{P}[C_n(\pmb{\mu},\pmb{\Theta}_n)^c \cap I_n(\pmb{\mu},\pmb{\Theta}_n)].
\end{align}

\vspace{-3mm}
In view of (\ref{eq:conn_ZeroLaw}), we obtain the zero-law for connectivity, i.e., that
\begin{equation} \nonumber
\lim_{n\to\infty} \mathbb{P}[\mathbb{G}(n;\pmb{\mu},\pmb{\Theta}_n) \text{ is connected}]= 0    \quad \text{ if } \quad c<1,
\end{equation}
immediately from the zero-law part of
Theorem \ref{theorem:isolated_nodes}, i.e., from that 
$\lim_{n\to\infty} \mathbb{P}[I_n(\pmb{\mu},\pmb{\Theta}_n)]=0$ if $c<1$.
It remains to establish the one-law for connectivity. In the remainder of this section, we assume that (\ref{scaling_condition_KG}) holds for some $c>1$. 
From Theorem \ref{theorem:isolated_nodes} and (\ref{eq:conn_OneLaw}), we see that the one-law for connectivity, i.e., that
\[
\lim_{n\to\infty} \mathbb{P}[\mathbb{G}(n;\pmb{\mu},\pmb{\Theta}_n) \text{ is connected}]= 1    \quad \text{ if } \quad c > 1,
\]
will follow if we show that
\begin{align}
\lim_{n \to \infty} \mathbb{P}[C_n(\pmb{\mu},\pmb{\Theta}_n)^c \cap I_n(\pmb{\mu},\pmb{\Theta}_n)] = 0.
\label{eq:conn_bounding}
\end{align}
Our approach will be to find a suitable upper bound for (\ref{eq:conn_bounding}) and prove that it goes to zero as $n$ goes to infinity with $c>1$.

We now work towards deriving an upper bound for (\ref{eq:conn_bounding}); then in {\color{\mycolor} Appendix~\ref{app:establishing}}
we will show that the bound goes to zero as $n$ gets large.
Define the event $E_n(\pmb{\mu},\pmb{\theta},\pmb{X})$ via
\begin{equation*}
E_n(\pmb{\mu},\pmb{\theta},\pmb{X}) :=\cup_{S \subseteq \mathcal{N} : |S| \geq 1} \left[ | \cup_{i \in S} \Sigma_i | \leq X_{|S|} \right]
\end{equation*}
where $\mathcal{N}=\{1,\ldots,n\}$ and $\pmb{X} = [X_1~ \cdots~ X_n]$ is an $n$-dimensional array of integers. Let
\begin{equation}
L_n := \min\left( \left \lfloor{\frac{P}{K_1}}\right \rfloor, \left \lfloor{\frac{n}{2}}\right \rfloor \right)
\label{eq:conn_newrange_Ln}
\end{equation}
and
\begin{equation}
X_\ell=
\begin{cases} 
      \hfill \left \lfloor{\beta \ell K_1}\right \rfloor    \hfill & \ell=1,\ldots,L_n \\
      \hfill \left \lfloor{\gamma P}\right \rfloor \hfill & \ell=L_n+1,\ldots,n \\
\end{cases}
\label{eq:conn_X}
\end{equation}
for some $\beta$ and $\gamma$ in $(0,\frac{1}{2})$ that will be specified later. In words, $E_n(\pmb{\mu},\pmb{\theta},\pmb{X})$ denotes the event that  there exists $\ell=1,\ldots, n$ such that the number of unique keys stored by at least one subset of $\ell$ sensors is less than $\left \lfloor{\beta \ell K_1}\right \rfloor \pmb{1}[\ell \leq L_n] + \left \lfloor{\gamma P}\right \rfloor \pmb{1}[\ell > L_n]$. Using a crude bound, we get
\begin{align}
&\mathbb{P}[C_n(\pmb{\mu},\pmb{\Theta}_n)^c \cap I_n(\pmb{\mu},\pmb{\Theta}_n)] \nonumber \\
& \leq \mathbb{P}[E_n(\pmb{\mu},\pmb{\theta}_n,\pmb{X}_n)] \nonumber \\
&\quad + \mathbb{P}[C_n(\pmb{\mu},\pmb{\Theta}_n)^c \cap I_n(\pmb{\mu},\pmb{\Theta}_n) \cap E_n(\pmb{\mu},\pmb{\theta}_n,\pmb{X}_n)^c] \label{eq:conn_2}
\end{align}
Thus, (\ref{eq:conn_bounding}) will be established by showing that
\begin{equation}
\lim_{n \to \infty} \mathbb{P}[E_n(\pmb{\mu},\pmb{\theta}_n,\pmb{X}_n)]=0, \label{eq:conn_3} \\
\end{equation}
and
\begin{equation}
\lim_{n \to \infty} \mathbb{P}[C_n(\pmb{\mu},\pmb{\Theta}_n)^c \cap I_n(\pmb{\mu},\pmb{\Theta}_n) \cap E_n(\pmb{\mu},\pmb{\theta}_n,\pmb{X}_n)^c] =0 \label{eq:conn_4}
\end{equation}

{
\prop
\label{prop:connectivity0}
Consider scalings $K_1,\ldots,K_r,P: \mathbb{N}_0 \rightarrow \mathbb{N}_0^{r+1}$ and $\alpha: \mathbb{N}_0 \rightarrow (0,1)$ such that (\ref{scaling_condition_KG}) 
holds for some $c>1$, (\ref{eq:conn_Pn2}) and (\ref{eq:conn_K1}) hold. Then, we have (\ref{eq:conn_3}) where $\pmb{X}_n$ is as specified in (\ref{eq:conn_X}), $\beta \in (0,\frac{1}{2})$ and $\gamma \in (0,\frac{1}{2})$ are selected such that
\begin{align}
&\max  \left( 2 \beta \sigma, \beta \left( \frac{e^2}{\sigma} \right)^{\frac{\beta}{1-2 \beta}} \right) < 1
\label{eq:conn_beta}
\\
&\max  \left( 2 \left( \sqrt{\gamma} \left( \frac{e}{\gamma} \right)^\gamma \right)^\sigma, \sqrt{\gamma} \left( \frac{e}{\gamma} \right)^\gamma \right) < 1
\label{eq:conn_gamma}
\end{align}
}

\begin{proof}
The proof is similar to \cite[Proposition~7.2]{Yagan/Inhomogeneous}. Results only require conditions (\ref{eq:conn_Pn2}) and $K_{1,n}= \omega(1)$ to hold. The latter condition is clearly established in Lemma~\ref{cor:new_corr}.
\end{proof}

The rest of the paper is devoted to establishing (\ref{eq:conn_4}) under the enforced assumptions on the scalings and with $\pmb{X}_n$ as specified in (\ref{eq:conn_X}), $\beta \in (0,\frac{1}{2})$ selected small enough such that (\ref{eq:conn_beta}) holds, and $\gamma \in (0,\frac{1}{2})$ selected small enough such that (\ref{eq:conn_gamma}) holds.
We denote by $\mathbb{G}(n,\pmb{\mu},\pmb{\Theta}_n)(S)$ a subgraph of $\mathbb{G}(n,\pmb{\mu},\pmb{\Theta}_n)$ whose vertices are restricted to the set $S$. Define the events

\vspace{-5mm}
\begin{align} \nonumber
 C_n(\pmb{\mu},\pmb{\Theta}_n,S) &:= [\mathbb{G}(n,\pmb{\mu},\pmb{\Theta}_n)(S) \text{ is connected}] 
\\
B_n(\pmb{\mu},\pmb{\Theta}_n,S) 
&:= [\mathbb{G}(n,\pmb{\mu},\pmb{\Theta}_n)(S) \text{ is isolated}]  \nonumber
\\
A_n(\pmb{\mu},\pmb{\Theta}_n,S) &:= C_n(\pmb{\mu},\pmb{\Theta}_n,S) \cap B_n(\pmb{\mu},\pmb{\Theta}_n,S)
\nonumber
\end{align}
In other words, $A_n(\pmb{\mu},\pmb{\Theta}_n,S)$ encodes the event that $\mathbb{G}(n,\pmb{\mu},\pmb{\Theta}_n)(S)$ is a \textit{component}, i.e., a connected subgraph that is isolated from the rest of the graph. The key observation is that a graph is {\em not} connected if and only if it has
a component on vertices $S$ with $1 \leq |S| \leq \left \lfloor{\frac{n}{2}}\right \rfloor$; note that if vertices $S$ form a component then so do vertices $\mathcal{N}-S$.
The event $I_n(\pmb{\mu},\pmb{\Theta}_n)$ eliminates the possibility of $\mathbb{G}(n,\pmb{\mu},\pmb{\Theta}_n)(S)$ containing a component of size one (i.e., an isolated node), whence we have
\begin{equation*}
C_n(\pmb{\mu},\pmb{\Theta}_n)^c \cap I_n(\pmb{\mu},\pmb{\Theta}_n) \subseteq \cup_{S \in \mathcal{N}:2 \leq |S| \leq \left \lfloor{\frac{n}{2}}\right \rfloor} A_n(\pmb{\mu},\pmb{\Theta}_n,S)
\end{equation*}
and the conclusion
\begin{equation}
\mathbb{P}[C_n(\pmb{\mu},\pmb{\Theta}_n)^c \cap I_n(\pmb{\mu},\pmb{\Theta}_n)] \leq \sum_{S \in \mathcal{N}:2 \leq |S| \leq \left \lfloor{\frac{n}{2}}\right \rfloor} \mathbb{P}[A_n(\pmb{\mu},\pmb{\Theta}_n,S)]
\nonumber
\end{equation}
follows. By exchangeability, we get
\begin{align}
&\mathbb{P}[C_n(\pmb{\mu},\pmb{\Theta}_n)^c \cap I_n(\pmb{\mu},\pmb{\Theta}_n) \cap E_n(\pmb{\mu},\pmb{\theta}_n,\pmb{X}_n)^c] 
\nonumber \\
&\leq \sum_{\ell=2}^{\left \lfloor{\frac{n}{2}}\right \rfloor} \left( \sum_{S \in \mathcal{N}_{n,\ell}} \mathbb{P}[A_n(\pmb{\mu},\pmb{\Theta}_n,S) \cap E_n(\pmb{\mu},\pmb{\theta}_n,\pmb{X}_n)^c] \right)
\nonumber \\
&= \sum_{\ell=2}^{\left \lfloor{\frac{n}{2}}\right \rfloor} \binom{n}{\ell} \mathbb{P}[A_{n,\ell}(\pmb{\mu},\pmb{\Theta}_n) \cap E_n(\pmb{\mu},\pmb{\theta}_n,\pmb{X}_n)^c]
\label{eq:key_bound1_connectivity_osy}
\end{align}
where $ \mathcal{N}_{n,\ell}$ denotes the collection of all subsets of $\{1,\ldots, n\}$ with exactly $\ell$ elements,
and $A_{n,\ell}(\pmb{\mu},\pmb{\Theta}_n)$ denotes the event that the set $\{1,\ldots, \ell \}$ of nodes form a component. 
As before we have $A_{n,\ell}(\pmb{\mu},\pmb{\Theta}_n) = C_\ell(\pmb{\mu},\pmb{\Theta}_n) \cap B_{n,\ell}(\pmb{\mu},\pmb{\Theta}_n) $,
where 
$C_\ell(\pmb{\mu},\pmb{\Theta}_n)$  denotes the event that $\{1,\ldots,\ell\}$ is connected and
$B_{n,\ell}(\pmb{\mu},\pmb{\Theta}_n)$ denotes the event that  $\{1,\ldots,\ell\}$ is isolated from the rest of the graph. 

Next, with $\ell=1,2,\ldots,n-1$, define $\nu_{\ell,j}(\alpha)$ by
\begin{equation}
\nu_{\ell,j}(\alpha):=\{i=1,2,\ldots,\ell : B_{ij}(\alpha)=1\}
\label{eq:isolated_v(alpha)}
\end{equation}
for each $j=\ell+1,\ldots,n$. Namely, $\nu_{\ell,j}(\alpha)$ is the set of nodes in $\{v_1,\ldots,v_{\ell}\}$ that are adjacent to node $v_j$ in the ER graph $\mathbb{H}(n;\alpha_n)$. For each $\ell=1,\ldots,n-1$, we have
\begin{equation*}
B_{n,\ell}(\pmb{\mu},\pmb{\Theta}_n) = \bigcap_{m=\ell+1}^{n} \left[\left( \cup_{i \in \nu_{\ell,m}(\alpha_n)} \Sigma_i \right) \cap \Sigma_m = \emptyset  \right].
\end{equation*}
We have
\begin{align}
& \bP{B_{n,\ell}(\pmb{\mu},\pmb{\Theta}_n) ~\big|~ \Sigma_1, \ldots, \Sigma_{\ell} 
 } 
\nonumber \\
&= \mathbb{E} \left[ \prod_{m=\ell+1}^n  \frac{\binom{P-| \cup_{i \in \nu_{\ell,m}(\alpha_n)} \Sigma_i |}{|\Sigma_m|}}{\binom {P}{|\Sigma_m|}} ~\bigg|~ \Sigma_1, \ldots, \Sigma_{\ell} \right]
\nonumber \\
&= \prod_{m=\ell+1}^n \mathbb{E} \left[ \frac{\binom{P-| \cup_{i \in \nu_{\ell,m}(\alpha_n)} \Sigma_i |}{|\Sigma_m|}}{\binom {P}{|\Sigma_m|}}
~\bigg|~ \Sigma_1, \ldots, \Sigma_{\ell} 
 \right]\nonumber \\
&=\mathbb{E} \left[ \frac{\binom{P-| \cup_{i \in \nu_{\ell}(\alpha_n)} \Sigma_i |}{\left|\Sigma\right|}}{\binom {P}{\left|\Sigma \right|}} ~\bigg|~ \Sigma_1, \ldots, \Sigma_{\ell} \right]^{n-\ell}
\label{eq:conn_6}
\end{align}
noting the fact that the collection of rvs $\{\nu_{\ell,m}, \Sigma_m: m=\ell+1,\ldots,n \}$ are mutually independent {\em and} identically distributed. Here,  $\nu_\ell(\alpha_n)$ denotes a generic rv distributed identically with $\nu_{\ell,m}(\alpha_n)$ for any $m=\ell+1,\ldots, n$.  Similarly, $\left|\Sigma\right|$ denotes a rv that takes the value $K_j$ with probability $\mu_j$. 

We will leverage the expression (\ref{eq:conn_6}) in (\ref{eq:key_bound1_connectivity_osy}) in the following manner.
Note that on the event $E_n(\pmb{\mu},\pmb{\theta}_n,\pmb{X}_n)^c$, we have
\begin{equation}
\left| \cup_{i \in \nu_\ell(\alpha_n)} \Sigma_i \right| \geq \left( X_{n,\nu_\ell(\alpha_n)} +1 \right) \pmb{1}[\left|\nu_\ell(\alpha_n)\right| >0]
\label{eq:conn_7}
\end{equation}
while the crude bound
\begin{equation}
\left| \cup_{i \in \nu_\ell(\alpha_n)} \Sigma_i \right| \geq K_{1,n} \pmb{1}[\left|\nu_\ell(\alpha_n)\right| >0]
\label{eq:conn_8}
\end{equation}
always holds. These bounds lead to
\begin{align}
& \bP{B_{n,\ell}(\pmb{\mu},\pmb{\Theta}_n)  \cap E_n(\pmb{\mu},\pmb{\theta}_n,\pmb{X}_n)^c~\big|~ \Sigma_1, \ldots, \Sigma_{\ell}}
 \nonumber \\
 & ~ \leq \mathbb{E} \left[ \frac{\binom{P-\max (K_{1,n},X_{n,\nu_\ell(\alpha_n)} +1)\pmb{1}[\left|\nu_\ell(\alpha_n)\right| >0]}{\left|\Sigma\right|}}{\binom {P}{\left|\Sigma\right|}} \right]^{n-\ell} 
\end{align}
Conditioning on $\Sigma_1, \ldots, \Sigma_{\ell}$ and $\{B_{ij}(\alpha_n), 1 \leq i < j \leq \ell \}$, we then get
\begin{align}
&\bP{A_{n,\ell}(\pmb{\mu},\pmb{\Theta}_n) \cap E_n(\pmb{\mu},\pmb{\theta}_n,\pmb{X}_n)^c} \nonumber \\
&=\bE{\1{C_\ell(\pmb{\mu},\pmb{\Theta}_n)} \1{B_{n,\ell}(\pmb{\mu},\pmb{\Theta}_n) \cap E_n(\pmb{\mu},\pmb{\theta}_n,\pmb{X}_n)^c}}\nonumber \\
& \leq \mathbb{E} \Bigg[ \pmb{1}[C_\ell(\pmb{\mu},\pmb{\Theta}_n)]  \cdot \nonumber \\
& \quad \cdot \mathbb{E} \left[ \left( \frac{\binom{P-\max (K_{1,n},X_{n,\nu_\ell(\alpha_n)} +1)\pmb{1}[\left|\nu_\ell(\alpha_n)\right| >0]}{\left|\Sigma\right|}}{\binom {P}{\left|\Sigma\right|}} \right) \Bigg]^{n-\ell} \right]\nonumber \\
&=\mathbb{P}[C_\ell(\pmb{\mu},\pmb{\Theta}_n)] \bE{\frac{\binom{P-\max (K_{1,n},X_{n,\nu_\ell(\alpha_n)} +1)\pmb{1}[\left|\nu_\ell(\alpha_n)\right| >0]}{\left|\Sigma\right|}}{\binom {P}{\left|\Sigma\right|}}}^{n-\ell}
\label{eq:conn_key_bound}
\end{align}
since  $C_\ell(\pmb{\mu},\pmb{\Theta}_n)$ is 
fully determined by $\Sigma_1, \ldots, \Sigma_{\ell}$ and $\{B_{ij}(\alpha_n), 1 \leq i < j \leq \ell \}$,
and $B_{n,\ell}(\pmb{\mu},\pmb{\Theta}_n)$ and $E_n(\pmb{\mu},\pmb{\theta}_n,\pmb{X}_n)$
are independent from $\{B_{ij}(\alpha_n), 1 \leq i, j \leq \ell \}$.

The next result establishes bounds for both terms  at (\ref{eq:conn_key_bound}).
{\lemma
\label{lemma:conn_key_lemma}
Consider a  distribution $\pmb{\mu}=(\mu_1,\mu_2,\ldots,\mu_r)$, integers 
$K_1  \leq \cdots  \leq K_r \leq P/2$, and  $\alpha \in (0,1)$. With $\pmb{X}_n$ as specified in (\ref{eq:conn_X}), $\beta \in (0,\frac{1}{2})$ and $\gamma \in (0,\frac{1}{2})$, we have
\begin{align}
\mathbb{P}[C_\ell(\pmb{\mu},\pmb{\Theta})] \leq \min \left\{ 1,\ell^{\ell-2} \left( \alpha p_{rr} \right)^{\ell-1} \right\}  
\label{eq:conn_key_bound_aux1}
\end{align}
and
\begin{align}
& \bE{\frac{\binom{P-\max (K_1,X_{n,\nu_\ell(\alpha)} +1)\pmb{1}[\left|\nu_\ell(\alpha)\right| >0]}{\left|\Sigma\right|}}{\binom {P}{\left|\Sigma\right|}}}
\leq  \min\bigg\{1-\alpha \lambda_1,
\nonumber
\\ \label{eq:conn_key_bound_aux2}
& ~ \min\{1-\mu_r + \mu_r e^{-\alpha p_{1r} \beta \ell}, e^{-\alpha p_{11} \beta \ell} \} + e^{- \gamma K_{1}} \pmb{1}[\ell>L_n]\bigg\}
\end{align}
}
The proof of Lemma~\ref{lemma:conn_key_lemma} is given in Appendix~\ref{app:key_lemma}.

Our proof of  (\ref{eq:conn_4}) will be completed (see (\ref{eq:key_bound1_connectivity_osy})) upon establishing 
\begin{equation}
\lim_{n \to \infty} \sum_{\ell=2}^{\left \lfloor{\frac{n}{2}}\right \rfloor} \binom{n}{\ell} \mathbb{P}[A_{n,\ell}(\pmb{\mu},\pmb{\Theta}_n) \cap E_n(\pmb{\mu},\pmb{\theta}_n,\pmb{X}_n)^c]=0
\label{eq:conn_key}
\end{equation}
by means of (\ref{eq:conn_key_bound}), (\ref{eq:conn_key_bound_aux1}), and (\ref{eq:conn_key_bound_aux2}).
These steps are taken in Appendix~\ref{app:establishing}. This establishes the one-law.

%

\ifCLASSOPTIONcaptionsoff
  \newpage
\fi

\vspace{-3mm}
\balance
\bibliographystyle{IEEEtran}
\bibliography{IEEEabrv,references}

\begin{thebibliography}{10}
\providecommand{\url}[1]{#1}
\csname url@samestyle\endcsname
\providecommand{\newblock}{\relax}
\providecommand{\bibinfo}[2]{#2}
\providecommand{\BIBentrySTDinterwordspacing}{\spaceskip=0pt\relax}
\providecommand{\BIBentryALTinterwordstretchfactor}{4}
\providecommand{\BIBentryALTinterwordspacing}{\spaceskip=\fontdimen2\font plus
\BIBentryALTinterwordstretchfactor\fontdimen3\font minus
  \fontdimen4\font\relax}
\providecommand{\BIBforeignlanguage}[2]{{%
\expandafter\ifx\csname l@#1\endcsname\relax
\typeout{** WARNING: IEEEtran.bst: No hyphenation pattern has been}%
\typeout{** loaded for the language `#1'. Using the pattern for}%
\typeout{** the default language instead.}%
\else
\language=\csname l@#1\endcsname
\fi
#2}}
\providecommand{\BIBdecl}{\relax}
\BIBdecl

\bibitem{Akyildiz_2002}
I.~Akyildiz, W.~Su, Y.~Sankarasubramaniam, and E.~Cayirci, ``A survey on sensor
  networks,'' \emph{IEEE Communications Magazine}, vol.~40, no.~8, pp.
  102--114, Aug 2002.

\bibitem{security_survey}
Y.~Wang, G.~Attebury, and B.~Ramamurthy, ``A survey of security issues in
  wireless sensor networks,'' \emph{IEEE Communications Surveys Tutorials},
  vol.~8, no.~2, pp. 2--23, Second 2006.

\bibitem{Gligor_2002}
L.~Eschenauer and V.~D. Gligor, ``A key-management scheme for distributed
  sensor networks,'' in \emph{Proc. of ACM CCS 2002}, pp. 41--47.

\bibitem{Haowen_2003}
H.~Chan, A.~Perrig, and D.~Song, ``Random key predistribution schemes for
  sensor networks,'' in \emph{Proc. of IEEE S\&P 2003}, May, pp. 197--213.

\bibitem{Raghavendra_2004}
C.~S. Raghavendra, K.~M. Sivalingam, and T.~Znati, Eds., \emph{Wireless Sensor
  Networks}.\hskip 1em plus 0.5em minus 0.4em\relax Norwell, MA, USA: Kluwer
  Academic Publishers, 2004.

\bibitem{camtepe_2005}
S.~A. \c{C}amtepe and B.~Yener, ``Key distribution mechanisms for wireless
  sensor networks: a survey,'' \emph{Rensselaer Polytechnic Institute, Troy,
  New York, Technical Report}, pp. 05--07, 2005.

\bibitem{Yarvis_2005}
M.~Yarvis, N.~Kushalnagar, H.~Singh, A.~Rangarajan, Y.~Liu, and S.~Singh,
  ``Exploiting heterogeneity in sensor networks,'' in \emph{Proc. of IEEE
  INFOCOM 2005}, March, pp. 878--890 vol. 2.

\bibitem{Yagan/Inhomogeneous}
O.~Ya\u{g}an, ``Zero-one laws for connectivity in inhomogeneous random key
  graphs,'' \emph{IEEE Transactions on Information Theory}, vol.~62, no.~8, pp.
  4559--4574, Aug 2016.

\bibitem{Gupta99}
P.~Gupta and P.~R. Kumar, ``Critical power for asymptotic connectivity in
  wireless networks,'' in \emph{Stochastic analysis, control, optimization and
  applications}.\hskip 1em plus 0.5em minus 0.4em\relax Springer, 1999, pp.
  547--566.

\bibitem{PenroseBook}
M.~D. Penrose, \emph{Random Geometric Graphs}.\hskip 1em plus 0.5em minus
  0.4em\relax {Oxford University Press}, Jul. 2003.

\bibitem{penrose2016connectivity}
M.~D. Penrose \emph{et~al.}, ``Connectivity of soft random geometric graphs,''
  \emph{The Annals of Applied Probability}, vol.~26, no.~2, pp. 986--1028,
  2016.

\bibitem{Yagan/EG_intersecting_ER}
O.~Ya\u{g}an, ``Performance of the {Eschenauer-Gligor} key distribution scheme
  under an {ON/OFF} channel,'' \emph{IEEE Transactions on Information Theory},
  vol.~58, no.~6, pp. 3821--3835, June 2012.

\bibitem{YaganPairwise}
O.~Ya\u{g}an and A.~Makowski, ``Modeling the pairwise key predistribution
  scheme in the presence of unreliable links,'' \emph{IEEE Transactions on
  Information Theory}, vol.~59, no.~3, pp. 1740--1760, March 2013.

\bibitem{yagan2012zero}
O.~Ya\u{g}an and A.~M. Makowski, ``Zero--one laws for connectivity in random
  key graphs,'' \emph{IEEE Transactions on Information Theory}, vol.~58, no.~5,
  pp. 2983--2999, 2012.

\bibitem{Zhao_2014}
J.~Zhao, O.~Ya\u{g}an, and V.~Gligor, ``On the strengths of connectivity and
  robustness in general random intersection graphs,'' in \emph{Proc. of IEEE
  CDC 2014}, 2014, pp. 3661--3668.

\bibitem{Rybarczyk}
M.~Bloznelis, J.~Jaworski, and K.~Rybarczyk, ``Component evolution in a secure
  wireless sensor network,'' \emph{Networks}, vol.~53, pp. 19--26, January
  2009.

\bibitem{Godehardt_2003}
E.~Godehardt and J.~Jaworski, ``Two models of random intersection graphs for
  classification,'' in \emph{Exploratory data analysis in empirical
  research}.\hskip 1em plus 0.5em minus 0.4em\relax Springer, 2003, pp. 67--81.

\bibitem{Bollobas}
B.~Bollob{\'a}s, \emph{Random graphs}.\hskip 1em plus 0.5em minus 0.4em\relax
  Cambridge university press, 2001, vol.~73.

\bibitem{DiPietroTissec}
R.~Di~Pietro, L.~V. Mancini, A.~Mei, A.~Panconesi, and J.~Radhakrishnan,
  ``Redoubtable sensor networks,'' \emph{Proc. of ACM TISSEC 2008}, vol.~11,
  no.~3, pp. 13:1--13:22, Mar. 2008.

\bibitem{MeiPanconesiRadhakrishnan2008}
A.~Mei, A.~Panconesi, and J.~Radhakrishnan, ``Unassailable sensor networks,''
  in \emph{Proc. of ACM SecureComm 2008}, pp. 1--10.

\bibitem{YM_ToN}
O.~Ya\u{g}an and A.~M. Makowski, ``Wireless sensor networks under the random
  pairwise key predistribution scheme: Can resiliency be achieved with small
  key rings?'' \emph{IEEE/ACM Transactions on Networking}, vol.~24, no.~6, pp.
  3383--3396, December 2016.

\bibitem{devroye2014connectivity}
L.~Devroye and N.~Fraiman, ``Connectivity of inhomogeneous random graphs,''
  \emph{Random Structures \& Algorithms}, vol.~45, no.~3, pp. 408--420, 2014.

\bibitem{Jun_2016}
J.~Zhao, ``Analyzing connectivity of heterogeneous secure sensor networks,''
  \emph{IEEE Transactions on Control of Network Systems}, vol.~PP, no.~99, pp.
  1--1, 2016.

\bibitem{zhaoconnectivity}
J.~Zhao, O.~Ya\u{g}an, and V.~Gligor, ``On connectivity and robustness in
  random intersection graphs,'' \emph{IEEE Transactions on Automatic Control},
  vol.~62, no.~5, pp. 2121--2136, May 2017.

\bibitem{Yagan_IP2013}
O.~Ya\u{g}an, D.~Qian, J.~Zhang, and D.~Cochran, ``Conjoining speeds up
  information diffusion in overlaying social-physical networks,'' \emph{IEEE
  Journal on Selected Areas in Communications}, vol.~31, no.~6, pp. 1038--1048,
  June 2013.

\bibitem{Bakshy_2012}
E.~Bakshy, I.~Rosenn, C.~Marlow, and L.~Adamic, ``The role of social networks
  in information diffusion,'' in \emph{Proc. of WWW 2012}.\hskip 1em plus 0.5em
  minus 0.4em\relax New York, NY, USA: ACM, pp. 519--528.

\bibitem{Yagan_Influence2012}
O.~Ya\ifmmode~\breve{g}\else \u{g}\fi{}an and V.~Gligor, ``Analysis of complex
  contagions in random multiplex networks,'' \emph{Phys. Rev. E}, vol.~86, p.
  036103, Sep 2012.

\bibitem{Kempe_2003}
D.~Kempe, J.~Kleinberg, and E.~Tardos, ``Maximizing the spread of influence
  through a social network,'' in \emph{Proc. of ACM SIGKDD 2003}.\hskip 1em
  plus 0.5em minus 0.4em\relax ACM, pp. 137--146.

\bibitem{Eubank_2004}
S.~Eubank, H.~Guclu, V.~Anil~Kumar, M.~V. Marathe, and e.~al,
  ``\BIBforeignlanguage{English}{Modeling disease outbreaks in realistic urban
  social networks},'' \emph{\BIBforeignlanguage{English}{Nature}}, vol. 429,
  no. 6988, pp. 180--4, May 13 2004.

\bibitem{Newman_2002}
M.~E.~J. Newman, ``Spread of epidemic disease on networks,'' \emph{Phys. Rev.
  E}, vol.~66, p. 016128, Jul 2002.

\bibitem{Barabasi_1999}
A.-L. Barab{\'a}si and R.~Albert, ``Emergence of scaling in random networks,''
  \emph{science}, vol. 286, no. 5439, pp. 509--512, 1999.

\bibitem{Watts_1998}
D.~J. Watts and S.~H. Strogatz, ``Collective dynamics
  of'small-world'networks,'' \emph{nature}, vol. 393, no. 6684, p. 440, 1998.

\bibitem{Amaral_2000}
L.~A.~N. Amaral, A.~Scala, M.~Barthelemy, and H.~E. Stanley, ``Classes of
  small-world networks,'' \emph{Proceedings of the National Academy of
  Sciences}, vol.~97, no.~21, pp. 11\,149--11\,152, 2000.

\bibitem{Yagan2017_SmallWorld}
O.~{Ya{\u g}an} and A.~M. {Makowski}, ``{Counting triangles, tunable clustering
  and the small-world property in random key graphs (Extended version)},''
  \emph{ArXiv e-prints}, Jan. 2017.

\bibitem{deijfen_kets_2009}
M.~Deijfen and W.~Kets, ``Random intersection graphs with tunable degree
  distribution and clustering,'' \emph{Probability in the Engineering and
  Informational Sciences}, vol.~23, no.~4, p. 661–674, 2009.

\bibitem{JansonLuczakRucinski}
S.~Janson, T.~{\L}uczak, and A.~Ruci{\'n}ski, \emph{Random Graphs. 2000}.\hskip
  1em plus 0.5em minus 0.4em\relax Wiley--Intersci. Ser. Discrete Math. Optim,
  2000.

\bibitem{jensen_99}
S.~T. Jensen, ``The laguerre-samuelson inequality with extensions and
  applications in statistics and matrix theory,'' Ph.D. dissertation,
  Department of Mathematics and Statistics, McGill University, 1999.

\bibitem{Yagan/PhD}
O.~Ya\u{g}an, ``Random graph modeling of key distribution scheme in wireless
  sensor networks,'' Ph.D. dissertation, University of Maryland, College Park
  (MD), 2011.

\bibitem{yagan_gradual}
O.~Ya\u{g}an and A.~M. Makowski, ``On the gradual deployment of random pairwise
  key distribution schemes,'' in \emph{Proc. of WiOpt 2011}, May 2011, pp.
  257--264.

\bibitem{Cayley}
G.~E. Martin, \emph{Counting: The art of enumerative combinatorics}.\hskip 1em
  plus 0.5em minus 0.4em\relax Springer Science \& Business Media, 2013.

\end{thebibliography}


\setcounter{equation}{0}
\renewcommand{\theequation}{\thesection.\arabic{equation}}

\newpage

\appendices

\section{Preliminaries}
\label{app:prelim}

{\lemma [{\cite[Lemma~4.2]{Yagan/Inhomogeneous}}]
\label{lemma:conn_asym_eq}
Consider any scaling $K_1,K_2,\ldots,K_r,P:\mathbb{N}_0 \rightarrow \mathbb{N}_0^{r+1}$. For any $i,j=1,\ldots,r$,
\begin{equation*}
\lim_{n \to \infty} p_{ij}(n)=0 \quad \text{if and only if} \quad \lim_{n \to \infty} \frac{K_{i,n} K_{j,n}}{P_n}=0
\end{equation*}
and we have the asymptotic equivalence
\begin{equation}
p_{ij}(n) \sim \frac{K_{i,n} K_{j,n}}{P_n}.
\label{eq:lemma9}
\end{equation}
}

\vspace{-6mm}
{\prop[{\cite[Proposition~4.4]{Yagan/Inhomogeneous}}]
For any set of positive integers $K_1,\ldots,K_r,P$ and any scalar $a \geq 1$, we have
\begin{equation}
\frac{\binom {P-\left \lceil{aK_i}\right \rceil }{K_j}}{\binom P{K_j}} \leq
\left(\frac{\binom {P-K_i}{K_j}}{\binom P{K_j}}\right)^a, \quad i,j=1,\ldots, r.
\label{eq:isolated:combinatorial_bound}
\end{equation}
}

\vspace{-3mm}
{\lemma
\label{lemma:isolated_1}
Consider a scaling $K_1,K_2,\ldots,K_r,P:\mathbb{N}_0 \rightarrow \mathbb{N}_0^{r+1}$ such that (\ref{scaling_condition_KG}) holds. We have
\begin{equation}
c_n \frac{\log n}{n \alpha_n} \leq p_{1r}(n) \leq \frac{c_n}{\mu_r} \frac{\log n}{n \alpha_n}
\label{eq:conn_p1r_asym}
\end{equation}
If in addition (\ref{eq:conn_K1}) holds, we have
\begin{equation}
p_{rr}(n)=o\left( \frac{(\log  n)^2}{n \alpha_n}\right).
\label{eq:conn_prr_asym}
\end{equation}
}
A proof of Lemma \ref{lemma:isolated_1} is given in Appendix \ref{sec:proof_bounds_p_ij}.

{
\lemma
\label{cor:new_corr}
Under (\ref{eq:conn_K1}), we have
\begin{equation}
\frac{K_{1,n}^2}{P_n}=\omega\left(\frac{1}{n \alpha_n}\right),
\label{eq:conn_K1_new_trick}
\end{equation}
and
\begin{equation}
K_{1,n}=\omega(1).
\label{eq:conn_K1_new_trick_2}
\end{equation}
}
\myproof
It is a simple matter to check that 
$
p_{11}(n) \leq \frac{K_{1,n}^2}{P_n-K_{1,n}}
$; see \cite[Proposition~7.1-7.2]{yagan2012zero} for a proof. 
In view of (\ref{scaling_condition_K}) this gives $p_{11}(n) \leq 2 \frac{K_{1,n}^2}{P_n}$.
Thus, we have
 \begin{equation} \nonumber
\frac{K_{1,n}^2}{P_n}=\Omega\left( p_{11}(n) \right)=\omega\left(\frac{1}{n \alpha_n}\right).
\end{equation}
From (\ref{eq:conn_Pn2}), (\ref{eq:conn_K1_new_trick}), and $\alpha_n \leq 1$, we readily obtain (\ref{eq:conn_K1_new_trick_2}).
\myendpf

\vspace{-3mm}
Other useful bounds that will be used throughout are
\begin{align}
& \binom{n}{\ell} \leq \left( \frac{en}{\ell} \right)^\ell, \quad \ell=1,\ldots,n, \quad n=1,2,\ldots
\label{eq:conn_bounds_1}
\\
& \sum_{\ell=2}^{\left \lfloor{\frac{n}{2}}\right \rfloor} \binom{n}{\ell} \leq 2^n
\label{eq:conn_bounds_2}
\end{align}

\section{A proof of Lemma \ref{lemma:isolated_1}}
\label{sec:proof_bounds_p_ij}

We know from (\ref{scaling_condition_KG_v2}) that
\begin{equation*}
\lambda_1(n)=\sum_{j=1}^r \mu_j p_{1j} =c_n  \frac{\log n}{\alpha_n n}.
\end{equation*}
Since $p_{1j}$ is monotone increasing in $j=1,\ldots,r$ by virtue of (\ref{eq:isolated_ordering_of_lambda}), we readily obtain the bounds
\begin{equation}
c_n \frac{\log n}{n \alpha_n} \leq p_{1r}(n) \leq \frac{c_n}{\mu_r} \frac{\log n}{n \alpha_n}
\label{eq:bound_on_p1r}
\end{equation}
which establishes (\ref{eq:conn_p1r_asym}).

In view of (\ref{eq:bound_on_p1r}) that 
implies $p_{1r}(n) = \Theta(\frac{\log n}{\alpha_n n})$, we will obtain (\ref{eq:conn_prr_asym})
if we show that $p_{rr}(n) = o( \log n )p_{1r}(n)$. Here this  will be established  by showing that
\begin{align}
p_{rr}(n) &\leq \max \left( 2,  \frac{8 c_n}{\mu_r} \frac{\log n}{w_n}\right) p_{1r}(n), \quad n=2,3, \ldots
\label{eq:to_show_p_rr_appendix}
\end{align}
 for some sequence $w_n$ such that $\lim_{n \to \infty} w_n = \infty$.
Fix $n=2,3,\ldots.$  We  have either
$p_{1r}(n) > \frac{1}{2}$,
or
$p_{1r}(n) \leq \frac{1}{2}$.
In the former case, it automatically holds that
\begin{equation}
p_{rr}(n) \leq 2 p_{1r}(n)
\label{eq:rig:1st_part}
\end{equation}
by virtue of the fact that $p_{rr}(n) \leq 1$.

Assume now that $p_{1r}(n) \leq \frac{1}{2}$.
We know from \cite[Lemmas~7.1-7.2]{yagan2012zero} that
\begin{equation}
1-e^{-\frac{K_{j,n}K_{r,n}}{P_n}} \leq p_{jr}(n) \leq \frac{K_{j,n}K_{r,n}}{P_n-K_{j,n}}, ~~ j=1,\ldots, r
\label{eq:rig_3}
\end{equation}
and it follows that
\begin{align}
\frac{K_{1,n}K_{r,n}}{P_n} \leq \log \left( \frac{1}{1-p_{1r}(n)} \right) \leq \log 2 < 1.
\label{eq:rig_4}
\end{align}
Using the fact that $1-e^{-x} \geq \frac{x}{2}$ with $x$ in $(0,1)$,
we then get
\begin{equation}
p_{1r}(n) \geq \frac{K_{1,n}K_{r,n}}{2P_n}.  
\label{eq:p_1r_lower_bound}
\end{equation}
In addition, using the upper bound in (\ref{eq:rig_3}) with $j=r$  gives
\[
p_{rr}(n) \leq  \frac{K_{r,n}^2}{P_n-K_{r,n}} \leq 2 \frac{K_{r,n}^2}{P_n}
\]
as we invoke (\ref{scaling_condition_K}). Combining the last two bounds we obtain
\begin{equation}
\frac{p_{rr}(n)}{p_{1r}(n)} \leq  4 \frac{K_{r,n}}{K_{1,n}}
\label{eq:bound_on_prr_by_p1r}
\end{equation}

In order to bound the term ${K_{r,n}}/{K_{1,n}}$, we recall from Lemma~\ref{cor:new_corr} that (\ref{eq:conn_K1}) implies (\ref{eq:conn_K1_new_trick}), i.e., that
$\frac{K_{1,n}^2}{P_n}=\frac{w_n}{n \alpha_n}$, for some sequence $w_n$ satisfying $\lim_{n \to \infty} w_n=\infty$. 
Using this together with (\ref{eq:p_1r_lower_bound}) and (\ref{eq:bound_on_p1r}) we then get
\begin{equation*}
 \frac{K_{r,n}}{K_{1,n}}=\frac{ \frac{K_{1,n}K_{r,n}}{P_n}}{\frac{K_{1,n}^2}{P_n}} \leq \frac{2p_{1r}(n)}{\frac{w_n}{n \alpha_n}} \leq \frac{2 \frac{c_n}{\mu_r} \frac{\log n}{n \alpha_n}}{\frac{w_n}{n \alpha_n}}=\frac{2 c_n}{\mu_r} \frac{\log n}{w_n}
\end{equation*}
Reporting this into (\ref{eq:bound_on_prr_by_p1r}) we get
\begin{align}
p_{rr}(n) &  \leq  \frac{8 c_n}{\mu_r} \frac{\log n}{w_n} p_{1r}(n).
\label{eq:rig:2nd_part}
\end{align}
Combining (\ref{eq:rig:1st_part}) and (\ref{eq:rig:2nd_part}), we readily obtain (\ref{eq:to_show_p_rr_appendix}).

\section{Proof of Proposition~\ref{prop:new_osy}}
\label{app:proof_of_prop_6.2}

Consider fixed $\pmb{\Theta}$. 
\begin{align*}
&\mathbb{E}\left[x_{n,1}(\pmb{\mu},\pmb{\Theta}) x_{n,2}(\pmb{\mu},\pmb{\Theta})\right] \nonumber \\
&=\mathbb{E}\left[\pmb{1}[v_1 \text{ is isolated }, v_2 \text{ is isolated} \cap t_1=1,t_2=1]\right]\nonumber \\
&=\mu_1^2 \mathbb{E}\left[\pmb{1}[v_1 \text{ is isolated }, v_2 \text{ is isolated}] \given[\Big]  t_1=1,t_2=1\right]\nonumber \\
&=\mu_1^2 \mathbb{E} \left[\pmb{1}[v_1 \nsim v_2]  \prod_{m=3}^n \pmb{1}[v_m \nsim v_1, v_m \nsim v_2] \given[\right] t_1=t_2=1\Bigg]
\end{align*}

Now we condition on $\Sigma_1 $ and $\Sigma_2$ and 
note that i) $\Sigma_1$ and $\Sigma_2$  determine $t_1$ and $t_2$;
and ii) the events $[v_1 \nsim v_2], \{[v_m \nsim v_1 \cap v_m \nsim v_2]\}_{m=3}^n$ are mutually independent given
$\Sigma_1$ and $\Sigma_2$.
Thus, we have
\begin{align} 
&\mathbb{E}[x_{n,1}(\pmb{\mu},\pmb{\Theta}) x_{n,2}(\pmb{\mu},\pmb{\Theta})]
 \nonumber \\
 &= \mu_1^2 \mathbb{E}\Bigg[\mathbb{P}\left[v_1 \nsim v_2 \given[\Big] \Sigma_1,\Sigma_2\right]  \cdot
  \label{eq:osy_new_zero_isol}
 \\
  & ~~~\prod_{m=3}^n \mathbb{P}\left[v_m \nsim v_1 \cap v_m \nsim v_2 \given[\Big] \Sigma_1,\Sigma_2\right]\given[\bigg] t_1=t_2=1 \bigg]
  \nonumber
 \end{align}
 
Define the $\{0,1\}$-valued rv $u(\pmb{\theta})$ by
\begin{equation}
u(\pmb{\theta}):=\pmb{1}[\Sigma_1 \cap \Sigma_2 \neq \emptyset].
\label{eq:isolated_u(theta)}
\end{equation}

Recalling (\ref{eq:isolated_v(alpha)}), (\ref{eq:osy_new_zero_isol}) gives
\begin{align*}
&\mathbb{E}[x_{n,1}(\pmb{\mu},\pmb{\Theta}) x_{n,2}(\pmb{\mu},\pmb{\Theta})] \nonumber \\
&= \mu_1^2 \mathbb{E}\Bigg[(1-\alpha)^{u(\pmb{\theta})}\prod_{m=3}^n \frac{\binom {P-\left|\cup_{i \in \nu_{2,m}(\alpha)} \Sigma_i\right|}{|\Sigma_m|}}{\binom P{|\Sigma_m|}} \given[\Bigg] t_1=t_2=1 \Bigg] \nonumber 
\end{align*}

Conditioned on $u(\pmb{\theta})=0$ and $v_1, v_2$ being class-1, we have
\begin{equation*} 
\left|\cup_{i \in \nu_{2,m}(\alpha)} \Sigma_i \right|=\left|\nu_{2,m}(\alpha)\right|K_1.
\end{equation*}
Also, we have
\[
\mathbb{P}[u(\pmb{\theta_n}) = 0 \given[] t_1=t_2=1] = 1-p_{11}.
\]
Thus, we get
\begin{align}
&\mathbb{E}[x_{n,1}(\pmb{\mu},\pmb{\Theta}) x_{n,2}(\pmb{\mu},\pmb{\Theta}) ~ \pmb{1}[u(\pmb{\theta})=0]] \nonumber \\
&=\mu_1^2 (1-p_{11}) \mathbb{E}\left[\prod_{m=3}^n \frac{\binom {P-| \nu_{2,m}(\alpha)  K_1 |}{|\Sigma_m|}}{\binom P{|\Sigma_m|}}\right]\nonumber \\
&=\mu_1^2 (1-p_{11}) \mathbb{E}\left[\frac{\binom {P-|\nu_{2,3}(\alpha)| K_1}{|\Sigma_3|}}{\binom P{|\Sigma_3|}}\right]^{n-2}
\nonumber \\
& = \mu _1^2 (1-p_{11})
 \left(\sum_{j=1}^r \mu_j \mathbb{E}\left[\frac{\binom {P-|\nu_{2,3}(\alpha)| K_1}{|\Sigma_3|}}{\binom P{|\Sigma_3|}} \given[\bigg] t_3 =j \right]\right)^{n-2}
 \nonumber \\
& = \mu _1^2 (1-p_{11})
 \left(\sum_{j=1}^r \mu_j \mathbb{E}\left[\frac{\binom {P-|\nu_{2,3}(\alpha)| K_1}{K_j}}{\binom P{K_j}} \right]\right)^{n-2}
 \nonumber \\ \nonumber
& \leq \mu_1^2 (1-p_{11}) \mathbb{E}\left[\sum_{j=1}^r \mu_j \left(\frac{\binom {P-K_1}{K_j}}{\binom P{K_j}}\right)^{|\nu_{2,3}(\alpha)|} \right]^{n-2} 
\end{align}
where we use  (\ref{eq:isolated:combinatorial_bound}) in the last step.

Now, let $Z(\pmb{\theta})$ denote a rv that takes the value
\begin{equation}
\frac{\binom {P-K_1}{K_j}}{\binom P{K_j}} \quad \text{ with probability } \quad \mu_j, \quad j=1,\ldots,r.
\label{eq:Z_teta_defn}
\end{equation}
In other words, $Z(\pmb{\theta}) = 1-p_{1j}$ with probability $\mu_j$ so that $\mathbb{E}[Z(\pmb{\theta})] = 1- \lambda_1$.
Then,
\begin{align}
&\mathbb{E}[x_{n,1}(\pmb{\mu},\pmb{\Theta}) x_{n,2}(\pmb{\mu},\pmb{\Theta}) \pmb{1}\left[u(\pmb{\theta})=0\right]] \nonumber \\
& \leq \mu_1^2 (1-p_{11}) \mathbb{E}\left[Z(\pmb{\theta})^{|\nu_{2,3}(\alpha)|} \right]^{n-2}
\label{eq:isolated_second_moment_6}
\end{align}

Under the independent on/off channel model, we have that $|\nu_{2,3}(\alpha)|$ is a Binomial rv, i.e., $|\nu_{2,3}(\alpha)| =_{st} \text{Bin}(2,\alpha)$. Hence,
\begin{align}
&\mu_1^2 (1-p_{11}) \mathbb{E}\left[Z(\pmb{\theta})^{|\nu_{2,3}(\alpha)|} \right]^{n-2} \label{eq:isolated_second_moment_7}
 \\
&=\mathbb{E}\left[ \sum_{i=0}^2 \binom{2}{i} \alpha^i (1-\alpha)^{2-i} Z(\pmb{\theta})^i\right]^{n-2}\nonumber \\
&=\mu_1^2 (1-p_{11}) \mathbb{E}\left[(1-\alpha)^2+2\alpha(1-\alpha)Z(\pmb{\theta})+\alpha^2 Z(\pmb{\theta})^2 \right]^{n-2}
\nonumber 
\end{align}

Conditioning on $u(\pmb{\theta})=1$ and $t_1=t_2=1$, we have
\begin{align*}
|\cup_{i \in \nu_{2,m}(\alpha)} \Sigma_i|  & =
\begin{cases} 
       0    \hfill & \text{ if $|\nu_{2,m}(\alpha)|=0$} \\
       K_1 \hfill & \text{ if $|\nu_{2,m}(\alpha)|=1$} \\
       2K_1-|\Sigma_1 \cap \Sigma_2| \hfill & \text{ if $|\nu_{2,m}(\alpha)|=2$} \\
\end{cases} \nonumber
\end{align*}
and by a crude bounding argument, we have
\begin{equation}
|\cup_{i \in \nu_{2,m}(\alpha)} \Sigma_i| \geq K_1 \pmb{1}[|\nu_{2,m}(\alpha)|>0]
\label{eq:isolated_second_moment_5}
\end{equation}
Using (\ref{eq:isolated_second_moment_5}) and recalling the analysis for $\mathbb{E}[x_{n,1}(\pmb{\mu},\pmb{\Theta}) x_{n,2}(\pmb{\mu},\pmb{\Theta}) \pmb{1}[u(\pmb{\theta})=0]]$, we obtain
\begin{align}
&\mathbb{E}[x_{n,1}(\pmb{\mu},\pmb{\Theta}) x_{n,2}(\pmb{\mu},\pmb{\Theta}) \pmb{1}[u(\pmb{\theta})=1]] \nonumber \\
&\leq \mu_1^2 p_{11} (1-\alpha) \mathbb{E}\left[Z(\pmb{\theta})^{\pmb{1}[|\nu_{2,3}(\alpha)|>0]} \right]^{n-2} \nonumber \\
&= \mu_1^2  p_{11} (1-\alpha) \mathbb{E}\left[(1-\alpha)^2+\left( 1-(1-\alpha)^2 \right) Z_n \right]^{n-2}
\label{eq:isolated_second_moment_8}
\end{align} 

Combining (\ref{eq:isolated_second_moment_6}), (\ref{eq:isolated_second_moment_7}), and (\ref{eq:isolated_second_moment_8}), we get
\begin{align}
&\mathbb{E}[x_{n,1}(\pmb{\mu},\pmb{\Theta}) x_{n,2}(\pmb{\mu},\pmb{\Theta})] \nonumber \\
&= \mathbb{E}[x_{n,1}(\pmb{\mu},\pmb{\Theta}) x_{n,2}(\pmb{\mu},\pmb{\Theta}) \left(\pmb{1}[u(\pmb{\theta})=0] + \pmb{1}[u(\pmb{\theta})=1]  \right)]\nonumber \\
& \leq \mu_1^2 (1-p_{11}) \Big((1-\alpha)^2+2\alpha(1-\alpha)\mathbb{E}[Z(\pmb{\theta})] \nonumber \\
&\quad +\alpha^2 \mathbb{E}\left[Z(\pmb{\theta})^2\right] \Big)^{n-2} + \mu_1^2 p_{11}(1-\alpha) \Big( (1-\alpha)^2 \nonumber \\
&\quad + \left(1-(1-\alpha)^2\right) \mathbb{E}[Z(\pmb{\theta})]\Big)^{n-2}
\label{eq:isolated_second_moment_9}
\end{align}

It is clear from (\ref{eq:int_isol_prob_osy}) and the definition of $Z(\pmb{\theta})$ that
\begin{align}
\mathbb{E}[x_{n,1}(\pmb{\mu},\pmb{\Theta})]&=\mu_1 \left( \sum_{j=1}^r \mu_j (1-\alpha p_{1j}) \right)^{n-1}\nonumber \\
&=\mu_1 \left((1-\alpha)+\alpha \mathbb{E}\left[Z(\pmb{\theta})\right] \right)^{n-1}
\label{eq:isolated_second_moment_10}
\end{align}

Combining (\ref{eq:isolated_second_moment_9}) and (\ref{eq:isolated_second_moment_10}), we get
\begin{align}
&\frac{\mathbb{E}[x_{n,1}(\pmb{\mu},\pmb{\Theta}) x_{n,2}(\pmb{\mu},\pmb{\Theta})] }{\mathbb{E}[x_{n,1}(\pmb{\theta})]^2} 
\nonumber \\
&   \leq  (1-p_{11})   \frac{\left((1-\alpha)^2+2\alpha(1-\alpha)\mathbb{E}[Z(\pmb{\theta})]+\alpha^2 \mathbb{E}[Z(\pmb{\theta})^2] \right)^{n-2}}{\left((1-\alpha)+\alpha \mathbb{E}\left[Z(\pmb{\theta})\right]  \right)^{2(n-1)}} \nonumber \\
& \ +  p_{11} \frac{\left( (1-\alpha)^2+ \left(1-(1-\alpha)^2\right) \mathbb{E}[Z(\pmb{\theta})]\right)^{n-2}}{\left((1-\alpha)+\alpha \mathbb{E}\left[Z(\pmb{\theta})\right]  \right)^{2(n-1)}} \nonumber \\
&:=A+B
\label{eq:defn_AB}
\end{align}
where we use the fact that $1-\alpha \leq 1$.

We now consider a scaling $\pmb{\Theta} : \mathbb{N}_0 \rightarrow \mathbb{N}_0^{r+1} \times (0,1)$
as stated in Proposition \ref{prop:new_osy} and
bound the terms $A$ and $B$ in turn. Our goal is to show that
\begin{equation}
\limsup_{n \to \infty} (A+B) \leq 1.
\label{eq:to_show_last_step_zero_law}
\end{equation}
First, we write $\mathbb{E}[Z(\pmb{\theta}_n)^2] = \mathbb{E}[Z(\pmb{\theta}_n)]^2+var[Z(\pmb{\theta}_n)]$,
where $var[Z(\pmb{\theta}_n)]$ can be bounded by the Popoviciu's inequality \cite[p. 9]{jensen_99} as follows
\begin{align*}
var[Z(\pmb{\theta}_n)] &\leq \frac{1}{4} \left( \max (Z(\pmb{\theta}_n))-\min(Z(\pmb{\theta}_n)) \right)^2\nonumber \\
&= \frac{1}{4} \left( \frac{\binom {P_n-K_{1,n}}{K_{1,n}}}{\binom{P_n}{K_{1,n}}} - \frac{\binom {P_n-K_{1,n}}{K_{r,n}}}{\binom{P_n}{K_{r,n}}} \right)^2 \nonumber \\
& \leq \frac{1}{4} \left( 1 - \frac{\binom {P_n-K_{1,n}}{K_{r,n}}}{\binom{P_n}{K_{r,n}}} \right)^2 = \frac{1}{4} \left(p_{1r}(n)\right)^2.
\end{align*}
Then, we get from the scaling condition (\ref{scaling_condition_KG_v2}) and (\ref{eq:bound_on_p1r}) that
\begin{align*}
\mathbb{E}[Z(\pmb{\theta}_n)^2] & \leq \mathbb{E}[Z(\pmb{\theta}_n)]^2 + \frac{1}{4} \left(\frac{c_n}{\mu_r} \frac{\log n}{n \alpha_n}\right)^2
\end{align*}

Reporting this into (\ref{eq:defn_AB}) we get  
\begin{align*}
&\hspace{-2mm}A \leq (1-p_{11}) \frac{\left( \left( (1-\alpha_n)+\alpha_n \mathbb{E}[Z(\pmb{\theta}_n)]\right)^2 + \left( \frac{c_n}{2 \mu_r} \frac{\log n}{n}\right)^2  \right)^{n-2}}{\left((1-\alpha_n)+\alpha_n \mathbb{E}\left[Z(\pmb{\theta}_n)\right]  \right)^{2(n-1)}}\nonumber \\
&=(1+o(1))(1-p_{11})   \left(1+\left(\frac{\frac{c_n}{2 \mu_r} \frac{\log n}{n} }{1-\alpha_n +\alpha_n \mathbb{E}\left[Z(\pmb{\theta}_n\right] }\right)^2\right)^{n-2}
\end{align*}
where we used the fact that
\begin{equation}
\left((1-\alpha_n)+\alpha_n \mathbb{E}[Z(\pmb{\theta}_n)]\right)^2 
= \left(1- \alpha_n \lambda_1(n) \right)^2=1-o(1) 
\label{eq:osy_second_moment_new_EZ}
\end{equation}
since 
$\alpha_n \lambda_1(n) = c_n \log n / n$. 
Finally, we have
\begin{align*}
&\left(1+\left(\frac{\frac{c_n}{2 \mu_r} \frac{\log n}{n} }{1-\alpha_n +\alpha_n \mathbb{E}\left[Z(\pmb{\theta}_n\right] }\right)^2\right)^{n-2}
\nonumber \\ 
& ~ \leq  \exp \left\{n \left(\frac{\frac{c_n}{2 \mu_r} \frac{\log n}{n} }{1-c_n \frac{\log n}{n} }\right)^2 \right\} = o(1)
\end{align*}
since $\lim_{n \to \infty} c_n = c >0$ and $\mu_r>0$. Thus, we obtain the bound 
\begin{equation}
A \leq \left(1-p_{11}\right) \left(1+o(1) \right).
\label{eq:bound_on_A}
\end{equation}

We now consider the second term  in (\ref{eq:defn_AB}). Recall (\ref{eq:osy_second_moment_new_EZ})
and that $\mathbb{E}\left[Z(\pmb{\theta}_n)\right] = 1- \lambda_1(n) = 1 -c_n \log n / n$.
 We have
\begin{align}
B
&=\frac{p_{11}}{\left(1-\alpha_n+\alpha_n \mathbb{E}\left[Z(\pmb{\theta}_n)\right]  \right)^2} \cdot \nonumber \\
& \quad \cdot \left(1+\frac{\alpha_n^2 \mathbb{E}[Z(\pmb{\theta}_n)] (1-\mathbb{E}[Z(\pmb{\theta}_n)])}{\left(1-\alpha_n+\alpha_n \mathbb{E}\left[Z(\pmb{\theta}_n)\right]  \right)^2}\right)^{n-2}\nonumber \\
&\leq p_{11} (1+o(1)) \exp \left\{n  \frac{\alpha_n^2 c_n \frac{\log n}{n \alpha_n} (1-c_n \frac{\log n}{n \alpha_n})}{(1-c_n \frac{\log n}{n})^2}\right\}
\nonumber \\
& \leq p_{11} (1+o(1))  \exp \left\{\frac{c_n \alpha_n \log n}{\left(1-c_n \frac{\log n}{n}\right)^2} \right\}
\label{eq:isolated:alpha-log-n}
\end{align}

We will now establish the desired result (\ref{eq:to_show_last_step_zero_law})
by using (\ref{eq:bound_on_A}) and (\ref{eq:isolated:alpha-log-n}). Our approach is based on the subsubsequence principle \cite[p. 12]{JansonLuczakRucinski} and considering the cases $\lim_{n \to \infty} \alpha_n \log n = 0 $ and $\lim_{n \to \infty} \alpha_n \log n \in (0, \infty]$ separately.

\paragraph{Assume that $\lim_{n \to \infty} \alpha_n \log n=0$}  From  (\ref{eq:isolated:alpha-log-n}) we get
$B \leq (1+o(1))  p_{11}$ and upon using (\ref{eq:bound_on_A}) we see that 
$A+B \leq  (1+o(1))$ establishing (\ref{eq:to_show_last_step_zero_law}) along subsequences with $\lim_{n \to \infty} \alpha_n \log n=0$.

\paragraph{Assume that $\lim_{n \to \infty} \alpha_n \log n \in (0,\infty]$}  Since $p_{1j}$ is monotonically increasing in $j=1,\ldots,r$ (see (\ref{eq:isolated_ordering_of_lambda})), we have
\begin{align*} 
\begin{split}
\lambda_1&=\sum_{j=1}^r \mu_j p_{1j}\geq p_{11} \sum_{j=1}^r \mu_j = p_{11}
\end{split}
\end{align*}
Thus, $p_{11} \leq \lambda_1(n) = c_n \log n /(\alpha_n n)$. Then, (\ref{eq:isolated:alpha-log-n}) gives
\begin{align*}
B &\leq \left(1+o(1) \right) \frac{c_n \log n}{\alpha_n n}  \exp\left\{\frac{c_n \alpha_n \log n}{\left(1-c_n {\log n}/{n}\right)^2}\right\} 
\nonumber \\
&=  \left(1+o(1) \right) \frac{c_n (\log n)^2}{\alpha_n \log n}    n^{-1+\frac{c_n \alpha_n }{\left(1-c_n {\log n}/{n}\right)^2}}
\nonumber \\
&=  o(1)
\end{align*}
since $\lim_{n \to \infty} \alpha_n \log n  > 0$ along this subsequence and
\[
\lim_{n \to \infty} -1+\frac{c_n \alpha_n }{\left(1- c_n\log n/{n}\right)^2} < 0 
\]
given that $\lim_{n \to \infty}c_n =c <1$.
From (\ref{eq:bound_on_A}) and the fact that $p_{11} \leq 1$, we have $A \leq 1+o(1)$, and
(\ref{eq:to_show_last_step_zero_law}) follows.

The two cases considered cover all the possibilities for the limit of $\alpha_n \log n$. By virtue of the subsubsequence principle \cite[p. 12]{JansonLuczakRucinski}, we get
(\ref{eq:to_show_last_step_zero_law}) without any condition on the sequence $\alpha_n \log n$; i.e., we obtain (\ref{eq:to_show_last_step_zero_law})
even when the sequence $\alpha_n \log n$ does not have a limit!

\section{Establishing Lemma~\ref{lemma:conn_key_lemma}}
\label{app:key_lemma}
The bounds given at Lemma \ref{lemma:conn_key_lemma}
are valid irrespective of how the parameters involved scale with $n$. 
Thus, we  consider fixed $\pmb{\Theta}$ with constraints given in the statement of Lemma \ref{lemma:conn_key_lemma}. 

We first establish (\ref{eq:conn_key_bound_aux2}) starting with the first bound. 
Recall that $\left|\nu_\ell(\alpha)\right| $ is a Binomial rv with $\ell$ trials and success probability $\alpha$.
Recall also the rv $Z(\pmb{\theta})$ defined at (\ref{eq:Z_teta_defn}).
Using a crude  bound and then (\ref{eq:isolated:combinatorial_bound})   we get 
\begin{align}
&\mathbb{E} \left[ \frac{\binom{P-\max (K_1,X_{n,\nu_\ell(\alpha)} +1)\pmb{1}[\left|\nu_\ell(\alpha)\right| >0]}{\left|\Sigma\right|}}{\binom {P}{\left|\Sigma\right|}} \right] \nonumber \\
&\leq \mathbb{E} \left[ \frac{\binom{P-K_1\pmb{1}[\left|\nu_\ell(\alpha)\right| >0]}{\left|\Sigma\right|}}{\binom {P}{\left|\Sigma\right|}} \right]\nonumber \\
& \leq \mathbb{E} \left[ Z(\pmb{\theta})^{\pmb{1}[\left|\nu_\ell(\alpha)\right| >0]} \right]\nonumber \\
&=\left(1-\alpha\right)^\ell + \left( 1- \left( 1-\alpha\right)^\ell\right) \mathbb{E}[Z(\pmb{\theta})]\nonumber \\
&\leq 1-\alpha+\alpha \bE{Z(\pmb{\theta})} = 1-\alpha \lambda_1(n).
\label{eq:conn_crude_bound_1}
\end{align}
upon noting that $ \bE{Z(\pmb{\theta})} = 1- \lambda_1 \leq 1$. 

Next, consider range $\ell=1,\ldots,L_n$, where we have
\begin{equation*}
\left(X_{n,\nu_\ell(\alpha)} +1\right) \pmb{1}[\left|\nu_\ell(\alpha)\right|>0] \geq \left \lceil{\beta \left|\nu_\ell(\alpha)\right| K_1}\right \rceil
\end{equation*}
Recalling (\ref{eq:isolated:combinatorial_bound}), we get
\begin{align}
&\mathbb{E} \left[ \frac{\binom{P-\max (K_1,X_{n,\nu_\ell(\alpha)} +1)\pmb{1}[\left|\nu_\ell(\alpha)\right| >0]}{\left|\Sigma\right|}}{\binom {P}{\left|\Sigma\right|}} \right] \nonumber \\
&\leq \mathbb{E} \left[ \frac{\binom{P-\beta \left|\nu_\ell(\alpha)\right| K_1}{\left|\Sigma\right|}}{\binom {P}{\left|\Sigma\right|}} \right]\nonumber \\
&= \mathbb{E}\left[ Z(\pmb{\theta})^{\beta \left|\nu_\ell(\alpha)\right|}\right]
\nonumber \\ 
 &= \mathbb{E}\left[ \sum_{j=0}^\ell \binom {\ell}{j} \alpha^j (1-\alpha)^{\ell-j} Z(\pmb{\theta})^{\beta j} \right]\nonumber \\
&=\mathbb{E} \left[ \left( 1-\alpha\left(1-Z(\pmb{\theta})^\beta \right) \right)^\ell \right]\nonumber \\
& \leq \mathbb{E} \left[ \left( 1-\alpha \beta \left(1-Z(\pmb{\theta})\right) \right)^\ell \right] \leq \mathbb{E} \left[ e^{-\alpha(1-Z(\pmb{\theta}))\beta \ell}\right]
\label{eq:conn_lemma_1}
\end{align}
using the fact 
that 
$1-Z(\pmb{\theta})^\beta \geq \beta(1-Z(\pmb{\theta}))$ 
with $Z(\pmb{\theta}) \leq 1$ and $0 \leq \beta \leq 1$; a proof is available at \cite[Lemma~5.2]{Yagan/EG_intersecting_ER}.
On the range $\ell=L_n+1,\ldots,\left \lfloor{\frac{n}{2}}\right \rfloor$, $\left|\nu_\ell(\alpha)\right|$ can be less than or greater than $L_n$. 
In the latter case, we have
\begin{equation*}
\max (K_1,X_{n,\nu_\ell(\alpha)} +1)\pmb{1}[\left|\nu_\ell(\alpha)\right| >0] \geq \left \lfloor{\gamma P}\right \rfloor +1
\end{equation*}
Using (\ref{eq:conn_lemma_1}) and the fact that (see \cite[Lemma~5.4.1]{Yagan/PhD} for a proof)
\begin{equation}\nonumber
{\binom{P-K_1}{K_2}}\bigg/{\binom{P}{K_2}} \leq e^{-\frac{K_2}{P}K_1}
\label{eq:isolated_BinoExpo}
\end{equation}
for $K_1+K_2 \leq P$,
we have
\begin{align}
&\mathbb{E} \left[ \frac{\binom{P-\max (K_1,X_{n,\nu_\ell(\alpha)} +1)\pmb{1}[\left|\nu_\ell(\alpha)\right| >0]}{\left|\Sigma\right|}}{\binom {P}{\left|\Sigma\right|}} \right] \nonumber \\
&\leq \mathbb{E} \left[ e^{-\alpha(1-Z(\pmb{\theta}))\beta \ell} \pmb{1}[\left|\nu_\ell(\alpha)\right| \leq L_n]  \right] \nonumber \\
& \quad + \mathbb{E} \left[ e^{-\frac{\left|\Sigma\right|}{P}(\left \lfloor{\gamma P}\right \rfloor +1)} \pmb{1}[\left|\nu_\ell(\alpha)\right| > L_n]  \right]
\nonumber \\ 
& \leq \mathbb{E}\left[e^{- \alpha (1-Z(\pmb{\theta})) \beta \ell}\right] + e^{- \gamma K_1} \pmb{1}[\ell>L_n]
\label{eq:conn_lemma}
\end{align}
by virtue of the fact that $ |\Sigma| \geq K_1$. 

Finally, we  get
(\ref{eq:conn_key_bound_aux2})
from  (\ref{eq:conn_crude_bound_1}) and (\ref{eq:conn_lemma}) by noting that
\begin{align}
\mathbb{E} \left[ e^{-\alpha(1-Z(\pmb{\theta}))\beta \ell} \right] 
&=\sum_{j=1}^r \mu_j  e^{-\alpha p_{1j} \beta \ell}
\leq (1-\mu_r)+\mu_re^{-\alpha p_{1r} \beta \ell}
\nonumber
\end{align}
and that
\begin{align}
\mathbb{E} \left[ e^{-\alpha(1-Z(\pmb{\theta}))\beta \ell} \right] &
= \sum_{j=1}^r \mu_j  e^{-\alpha p_{1j} \beta \ell}
\leq e^{-\alpha p_{11} \beta \ell}
\label{eq:conn_crude_bound_2}
\end{align}
 The last step used the fact that $p_{ij}$ is monotone increasing in both $i$ and $j$.
 
Next, we establish (\ref{eq:conn_key_bound_aux1}). This is a version of a fairly standard bound derived previously for various other random graph models including ER graphs \cite{Bollobas}, random key graphs  \cite{yagan2012zero}, and random $K$-out graphs \cite{yagan_gradual,YaganPairwise}.
The proof is very similar to that of \cite[Proposition~9.1]{Yagan/Inhomogeneous} and \cite[Lemma~10.2]{Yagan/EG_intersecting_ER}. We give it below for completeness.

Let $\mathbb{G}_\ell(n;\pmb{\mu},\pmb{\Theta})$ denote the subgraph of $\mathbb{G}(n;\pmb{\mu},\pmb{\Theta})$ induced on the vertices $\{ v_1,\ldots,v_\ell \}$. $\mathbb{G}_\ell(n;\pmb{\mu},\pmb{\Theta})$ is connected if and only if it contains a spanning tree; i.e., we have
\begin{equation} \nonumber
C_\ell(\pmb{\mu},\pmb{\Theta}) = \cup_{T \in \mathcal{T}_\ell} \left[ T \subseteq \mathbb{G}_\ell(n;\pmb{\mu},\pmb{\Theta}) \right]
\end{equation}
where $\mathcal{T}_\ell$ denotes the collection of all spanning trees on the vertices $\{ v_1,\ldots,v_\ell \}$. 
Thus,
\begin{equation} 
\mathbb{P}[C_\ell(\pmb{\mu},\pmb{\Theta})] \leq \sum_{T \in \mathcal{T}_\ell} \mathbb{P}\left[ T \subseteq \mathbb{G}_\ell(n;\pmb{\mu},\pmb{\Theta}) \right].
\label{eq:app_b_new1}
\end{equation}
Given that $K_1 \leq K_2 \leq \ldots \leq K_r$, the probability of $T$ being contained in $\mathbb{G}_\ell(n;\pmb{\mu},\pmb{\Theta})$ is maximized when all nodes receive the largest possible number $K_r$ of keys.
Thus, for any $T \in \mathcal{T}$ and  distribution $\pmb{\mu}$ we have
\begin{align} 
\mathbb{P}\left[ T \subseteq \mathbb{G}_\ell(n;\pmb{\mu},\pmb{\Theta}) \right] 
& \leq \mathbb{P}\left[ T \subseteq \mathbb{G}_\ell(n;\pmb{\mu}=\{0,0,\ldots,1\},\pmb{\Theta}) \right] \nonumber
\nonumber \\
& = (\alpha p_{rr})^{\ell-1}
\label{eq:app_b_new2}
\end{align}
where the last equality follows from the facts that i) a tree on $\ell$ vertices contain $\ell-1$ edges, 
and ii) since all nodes have the same key ring size, edges in  
$\mathbb{G}_{\ell} (n;\pmb{\mu}=\{0,0,\ldots,1\},\pmb{\Theta})$ are {\em pairwise} independent; see
\cite[Lemma~9.1]{yagan2012zero} and \cite[Eq. 64]{Yagan/EG_intersecting_ER}.
We obtain (\ref{eq:conn_key_bound_aux1}) upon using (\ref{eq:app_b_new2}) in (\ref{eq:app_b_new1}) and noting
by Cayley's formula \cite{Cayley} that
there are $\ell^{\ell-2}$ trees on $\ell$ vertices, i.e.,  $|\mathcal{T}_{\ell}| = \ell^{\ell-2}$.

\section{Establishing (\ref{eq:conn_key})}
\label{app:establishing}
We will establish (\ref{eq:conn_key}) in several steps with each step focusing on a specific range of the summation over $\ell$. Throughout, we consider a scalings $K_1,\ldots,K_r,P: \mathbb{N}_0 \rightarrow \mathbb{N}_0^{r+1}$ and $\alpha: \mathbb{N}_0 \rightarrow (0,1)$ such that (\ref{scaling_condition_KG}) 
holds with $c>1$, (\ref{eq:conn_K1}), and (\ref{eq:conn_Pn2}) hold.
\subsubsection{The case where $2 \leq \ell \leq R$}
This range considers fixed values of $\ell$. Pick an integer $R$ to be specified later at (\ref{eq:choosing_R}). Use (\ref{scaling_condition_KG}), (\ref{eq:conn_prr_asym}), (\ref{eq:conn_bounds_1}), (\ref{eq:conn_key_bound}), (\ref{eq:conn_key_bound_aux1}), and the first bound in (\ref{eq:conn_key_bound_aux2}) to get
\begin{align*}
&\sum_{\ell=2}^{R} \binom{n}{\ell} \mathbb{P}[A_{n,\ell}(\pmb{\mu},\pmb{\Theta}_n) \cap E_n(\pmb{\mu},\pmb{\theta}_n,\pmb{X}_n)^c] \nonumber \\
&\leq \sum_{\ell=2}^{R} \left( \frac{en}{\ell}\right)^\ell  \ell^{\ell-2} \left( \alpha_n p_{rr}(n)\right)^{\ell-1} \left( 1- \alpha_n \lambda_1(n) \right)^{n-\ell}\nonumber \\
& \leq \sum_{\ell=2}^{R} \left( en \right)^\ell  \left( \frac{(\log n)^2}{n}\right)^{\ell-1} \left( 1-c_n \frac{\log n}{n}\right)^{n-\ell}\nonumber \\
&\leq \sum_{\ell=2}^{R}  n \left( e(\log n)^2\right)^{\ell} e^{-c_n \log n \frac{n-\ell}{n}}\nonumber \\
&=\sum_{\ell=2}^{R} \left( e(\log n)^2\right)^{\ell}  n^{1-c_n\frac{n-\ell}{n}} \nonumber
\end{align*}
With $c>1$, we have $\lim_{n \to \infty} \left(1-c_n\frac{n-\ell}{n}\right)=1-c<0$.
Thus, 
for each $\ell=2,3, \ldots$, we have
\begin{equation} \nonumber
\left( e(\log n)^2\right)^{\ell-1}  n^{1-c_n\frac{n-\ell}{n}}=o(1),
\end{equation}
whence we get
\begin{equation} \nonumber
\lim_{n \to \infty} \sum_{\ell=2}^{R} \binom{n}{\ell} \mathbb{P}[A_{n,\ell}(\pmb{\mu},\pmb{\Theta}_n) \cap E_n(\pmb{\mu},\pmb{\theta}_n,\pmb{X}_n)^c] =0.
\end{equation}

\subsubsection{The case where $R+1 \leq \ell \leq \min\{L_n,  \lfloor{\frac{\mu_r n}{\beta c_n \log n}} \rfloor \}$}
\label{subsection:geometric_series}
Our goal in this and the next subsubsection is to cover the range $R+1 \leq \ell \leq  \lfloor{\frac{\mu_r n}{\beta c_n \log n}} \rfloor $.
Since the bound given at (\ref{eq:conn_key_bound_aux2}) takes a different form when $\ell > L_n$, we first consider the 
range $R+1 \leq \ell \leq \min\{L_n,  \lfloor{\frac{\mu_r n}{\beta c_n \log n}} \rfloor \}$.
Using (\ref{eq:conn_prr_asym}), (\ref{eq:conn_bounds_1}), (\ref{eq:conn_key_bound}), 
(\ref{eq:conn_key_bound_aux1}), and the second bound in (\ref{eq:conn_key_bound_aux2}) we get
\begin{align}
&\sum_{\ell=R+1}^{\min\{L_n, \lfloor{\frac{\mu_r n}{\beta c_n \log n}} \rfloor\}} \binom{n}{\ell} \mathbb{P}[A_{n,\ell}(\pmb{\mu},\pmb{\Theta}_n) \cap E_n(\pmb{\mu},\pmb{\theta}_n,\pmb{X}_n)^c] \nonumber \\
& \leq \sum_{\ell=R+1}^{\min\{L_n, \lfloor{\frac{\mu_r n}{\beta c_n \log n}} \rfloor \}}  \left( \frac{en}{\ell}\right)^\ell  \ell^{\ell-2} 
\left( \frac{(\log n)^2}{n}\right)^{\ell-1} \cdot
 \label{eq:conn_11}
\\ \nonumber
& \qquad \qquad \qquad \qquad \cdot \Bigg( 1 -\mu_r \left( 1-e^{- \alpha_n \beta \ell p_{1r}(n)} \right) \Bigg)^{n-\ell}
\end{align}
From the upper bound in (\ref{eq:conn_p1r_asym}) and $\ell \leq \frac{\mu_r n}{ \beta c_n \log n}$, we have 
\begin{align*}
\alpha_n \beta \ell p_{1r}(n)  &\leq \alpha_n \beta \frac{\mu_r n}{ \beta c_n \log n} \frac{c_n}{\mu_r} \frac{\log n}{n \alpha_n}  =  1.
\end{align*}
Using the fact that
$
1-e^{-x} \geq \frac{x}{2}$ for all $0 \leq x \leq 1$,
we  get
\begin{equation}
1-\mu_r \hspace{-.7mm}\left( 1-e^{- \alpha_n \beta \ell p_{1r}(n)} \right) \hspace{-.7mm}\leq 1  -  \frac{\mu_r \alpha_n \beta \ell p_{1r}(n)}{2} \leq \hspace{-.5mm} e^{- \beta \ell c_n \mu_r \frac{\log n}{2 n} }
\label{eq:conn_12}
\end{equation}
using the lower bound in (\ref{eq:conn_p1r_asym}).
Reporting this last bound in to (\ref{eq:conn_11}) and noting that 
\begin{equation}
n-\ell \geq \frac{n}{2}, \qquad \ell = 2, 3, \ldots, \left \lfloor \frac{n}{2} \right \rfloor,
\label{eq:bounding_n_minus_l}
\end{equation}
we get
\begin{align}
&\sum_{\ell=R+1}^{\min\{L_n,\left \lfloor{\frac{\mu_r n}{\beta c_n \log n}}\right \rfloor\}} \binom{n}{\ell} \mathbb{P}[A_{n,\ell}(\pmb{\mu},\pmb{\Theta}_n) \cap E_n(\pmb{\mu},\pmb{\theta}_n,\pmb{X}_n)^c] \nonumber \\
&\leq \sum_{\ell=R+1}^{\min\{L_n, \lfloor{\frac{\mu_r n}{\beta c_n \log n}} \rfloor\}}
 n \left( e(\log n)^2\right)^{\ell} 
e^{- \beta \ell c_n \mu_r\frac{\log n}{2n} \frac{n}{2}}\nonumber \\
&\leq n \sum_{\ell=R+1}^{\min\{L_n, \lfloor{\frac{\mu_r n}{\beta c_n \log n}} \rfloor\}} \left( e \left( \log n\right)^2 e^{- \beta c_n  \frac{\mu_r}{4} \log n} \right)^\ell\nonumber \\
&\leq n \sum_{\ell=R+1}^{\infty} \left( e \left( \log n\right)^2 e^{- \beta c_n  \frac{\mu_r}{4} \log n} \right)^\ell
\label{eq:conn_13}
\end{align}

Given that 
 $\beta, \mu_r>0$ and $\lim_{n \to \infty} c_n = c > 0$
we clearly have
\begin{equation}
 e \left( \log n\right)^2 e^{- \beta c_n \log n \frac{\mu_r}{4}} = o(1).
 \label{eq:summable_sequence}
 \end{equation}
Thus, the geometric series in (\ref{eq:conn_13}) is summable, and we have
\begin{align*}
&\sum_{\ell=R+1}^{\min\{L_n, \lfloor{\frac{\mu_r n}{\beta c_n \log n}} \rfloor\}} \binom{n}{\ell} \mathbb{P}[A_{n,\ell}(\pmb{\mu},\pmb{\Theta}_n) \cap E_n(\pmb{\mu},\pmb{\theta}_n,\pmb{X}_n)^c] \nonumber \\
&\leq \left(1+o(1)\right) n \left( e \left( \log n\right)^2 e^{- \beta c_n \log n \frac{\mu_r}{4}} \right)^{R+1}\nonumber \\
&=\left(1+o(1)\right) n^{1-(R+1){\beta c_n}\frac{\mu_r}{4}} \left( e (\log n)^2 \right)^{R+1} \nonumber \\
&= o(1)
\end{align*}
for any positive integer $R$ with
\begin{equation} 
R>\frac{8}{\beta c \mu_r}. 
\label{eq:choosing_R}
\end{equation}
This choice is permissible given that $c, \beta, \mu_r > 0$. 

\subsubsection{The case where $ \min \{ \lfloor{\frac{\mu_r n}{ \beta c_n \log n}} \rfloor,\max (R,L_n) \} < \ell \leq  \lfloor{\frac{\mu_r n}{ \beta c_n \log n}} \rfloor$}

Clearly, this range becomes obsolete if $\max (R,L_n) \geq  \lfloor{\frac{\mu_r n}{ \beta c_n \log n}} \rfloor$. Thus, it suffices to consider the subsequences for which the range $\max (R,L_n)+1 \leq \ell \leq  \lfloor{\frac{\mu_r n}{\beta c_n \log n}} \rfloor$ is non-empty. There, we use  (\ref{eq:conn_prr_asym}), (\ref{eq:conn_bounds_1}), (\ref{eq:conn_key_bound}), 
(\ref{eq:conn_key_bound_aux1}), and the second bound in (\ref{eq:conn_key_bound_aux2}) to get
\begin{align}
&\sum_{\ell=\max (R,L_n)+1}^{\left \lfloor{\frac{\mu_r n}{\beta c_n \log n}}\right \rfloor} \binom{n}{\ell}  \mathbb{P}[A_{n,\ell}(\pmb{\mu},\pmb{\Theta}_n) \cap E_n(\pmb{\mu},\pmb{\theta}_n,\pmb{X}_n)^c] \label{eq:conn_newrange_2}
 \\
&\leq \sum_{\ell=\max (R,L_n)+1}^{\left \lfloor{\frac{\mu_r n}{ \beta c_n \log n}}\right \rfloor} \left( \frac{en}{\ell} \right)^\ell \ell^{\ell-2} \left( \frac{\left(\log n \right)^2}{n} \right)^{\ell-1} . \nonumber \\
& \left( 1-\mu_r \left(1-e^{-\beta \ell \alpha_n p_{1r}(n)} \right)+e^{-\gamma K_{1,n}} \right)^{\frac{n}{2}}.
\nonumber \\
& \leq \hspace{-2mm} \sum_{\ell=\max (R,L_n)+1}^{\left \lfloor{\frac{\mu_r n}{2 \beta c \log n}}\right \rfloor} \hspace{-3mm} n \left(e \left(\log n \right)^2 \right)^\ell \left( e^{- \beta \ell c_n \mu_r \frac{\log n}{2n} } + e^{-\gamma K_{1,n}}\right)^{\frac{n}{2}}
\nonumber
\end{align}
where in the last step we used
 (\ref{eq:conn_12}) in view of $\ell \leq \frac{\mu_r n}{\beta c_n \log n}$.
 
Next, we write
\begin{align}
&e^{- \beta \ell c_n \mu_r \frac{\log n}{2n} } + e^{-\gamma K_{1,n}} \nonumber \\
&=  e^{- \beta \ell c_n \mu_r \frac{\log n}{2n} } \left(1+e^{-\gamma K_{1,n} + \beta \ell c_n \mu_r \frac{\log n}{2n}} \right) \nonumber \\
& \leq \exp \left\{- \beta \ell c_n \mu_r \frac{\log n}{2n} + e^{-\gamma K_{1,n} + \beta \ell c_n \mu_r \frac{\log n}{2n}} \right\} \nonumber \\
& \leq \exp \left\{- \beta \ell c_n  \mu_r \frac{\log n}{2n}\left(1 - \frac{e^{-\gamma K_{1,n} + \frac{\mu_r^2 }{2}}}{\beta \ell c_n \mu_r \frac{\log n}{2n}} \right) \right\} \label{eq:intermediary_range_new}
\end{align}
where the last inequality is obtained from $\ell \leq {\frac{\mu_r n}{ \beta c_n \log n}}$. 
Using the fact that $\ell > L_n=\min\{  \lfloor{\frac{P_n}{K_{1,n}}} \rfloor,  \lfloor{\frac{n}{2}} \rfloor \}$
 and (\ref{eq:conn_Pn2}) we have 
\begin{align}
\frac{e^{-\gamma K_{1,n}}}{\beta \ell c_n \mu_r \frac{\log n}{2n}}  & \leq 
\max\left\{  \frac{K_{1,n}}{P_n},  \frac{2}{n} \right \}{2n}
\frac{e^{-\gamma K_{1,n}}}{\beta  c_n \mu_r \log n} \nonumber \\
&\leq \max \left \{ \frac{2 K_{1,n} e^{-\gamma K_{1,n}}}{\beta  c_n \mu_r \sigma \log n},\frac{4 e^{-\gamma K_{1,n}}}{\beta  c_n \mu_r \log n} \right\} \nonumber \\
&=o(1) \nonumber
\end{align}
by virtue of (\ref{eq:conn_K1_new_trick_2}) and  the facts that $\beta, \mu_r, \sigma, c_n >0$. Reporting this
into (\ref{eq:intermediary_range_new}), we see that for
for any $\epsilon>0$, there exists a finite integer $n^*(\epsilon)$ such that 
\begin{equation}
\left( e^{- \beta \ell c_n \mu_r \frac{\log n}{2n}} + e^{-\gamma K_{1,n}}\right) \leq e^{- \beta \ell c_n \mu_r \frac{\log n}{2n} (1-\epsilon)} 
\label{eq:conn_newrange_3}
\end{equation}
for all $n \geq n^*(\epsilon)$.
Using (\ref{eq:conn_newrange_3}) in (\ref{eq:conn_newrange_2}), we get
\begin{align}
&\sum_{\ell=\max(R,L_n)+1}^{\left \lfloor{\frac{\mu_r n}{ \beta c_n \log n}}\right \rfloor} \binom{n}{\ell}  \mathbb{P}[A_{n,\ell}(\pmb{\mu},\pmb{\Theta}_n) \cap E_n(\pmb{\mu},\pmb{\theta}_n,\pmb{X}_n)^c] \nonumber \\
&\leq n \sum_{\ell=\max(R,L_n)+1}^{\left \lfloor{\frac{\mu_r n}{ \beta c_n \log n}}\right \rfloor} \left( e \left(\log n\right)^2 e^{- \beta  c_n \mu_r \frac{\log n}{2n} (1-\epsilon) \frac{n}{2}}  \right)^\ell \nonumber \\
&\leq n \sum_{\ell=\max (R,L_n)+1}^{\infty} \left( e \left(\log n\right)^2 e^{- \beta  c_n \mu_r \frac{\log n}{4} (1-\epsilon)}  \right)^\ell
\label{eq:conn_newrange_4}
\end{align}
Similar to (\ref{eq:summable_sequence}), we have
$e \left(\log n\right)^2 e^{- \beta  c_n \mu_r \frac{\log n}{4} (1-\epsilon)} = o(1)$ so that the sum 
in (\ref{eq:conn_newrange_4}) converges. 
Following a similar approach to that in Section~\ref{subsection:geometric_series}, we then see that 
\begin{equation*}
\lim_{n \to \infty} \sum_{\ell=\max (R,L_n)+1}^{\left \lfloor{\frac{\mu_r n}{2 \beta c \log n}}\right \rfloor} \binom{n}{\ell} \mathbb{P}[A_{n,\ell}(\pmb{\mu},\pmb{\Theta}_n) \cap E_n(\pmb{\mu},\pmb{\theta}_n,\pmb{X}_n)^c] = 0   
\end{equation*}
with $R$ selected according to (\ref{eq:choosing_R}) and $\epsilon<1/2$.

\subsubsection{The case where $\lfloor{\frac{\mu_r n}{\beta c_n \log n}} \rfloor+1 \leq \ell \leq \left \lfloor{\nu n}\right \rfloor $} 
\label{subsec:range_4}
We  consider $\lfloor{\frac{\mu_r n}{ \beta c_n \log n}} \rfloor+1 \leq \ell \leq \left \lfloor{\nu n}\right \rfloor$ for some $\nu \in (0,\frac{1}{2})$ to be specified later. Recall (\ref{eq:conn_p1r_asym}), (\ref{eq:conn_bounds_1}), (\ref{eq:conn_key_bound}), the first bound in (\ref{eq:conn_key_bound_aux1}), and the second bound in (\ref{eq:conn_key_bound_aux2}).
Noting that $\binom{n}{\ell}$ is monotone increasing in $\ell$ when $0 \leq \ell \leq \left \lfloor{\frac{n}{2}}\right \rfloor$  
 and using (\ref{eq:bounding_n_minus_l}) we get
\begin{align}
&\sum_{\ell= \lfloor{\frac{\mu_r n}{\beta c_n \log n}} \rfloor +1}^{ \left \lfloor{\nu n}\right \rfloor} \binom{n}{\ell} \mathbb{P}[A_{n,\ell}(\pmb{\mu},\pmb{\Theta}_n) \cap E_n(\pmb{\mu},\pmb{\theta}_n,\pmb{X}_n)^c] \nonumber \\
&\leq  \hspace{-3mm}\sum_{\ell= \lfloor{\frac{\mu_r n}{ \beta c_n \log n}} \rfloor+1}^{\left \lfloor{\nu n} \right \rfloor} \hspace{-2mm}\binom{n}{\left \lfloor{\nu n}\right \rfloor}  \hspace{-1.4mm}\left( 1-\mu_r+\mu_re^{- \alpha_n \beta \ell p_{1r}(n)} + e^{- \gamma K_{1,n}}\right)^{\frac{n}{2}}\nonumber \\
& \leq \hspace{-2mm} \sum_{\ell= \lfloor{\frac{\mu_r n}{ \beta c_n \log n}} \rfloor+1}^{\left \lfloor{\nu n}\right \rfloor} \hspace{-1mm} \left( \frac{e}{\nu} \right)^{\nu n} \bigg( 1-\mu_r +\mu_re^{- \alpha_n \beta  \frac{\mu_r n}{ \beta c_n \log n}  \frac{c_n \log n}{ n \alpha_n}}
\nonumber \\ 
 & ~~~~~~~~~ \qquad \qquad \qquad ~~~~  +\hspace{-.6mm}  e^{- \gamma K_{1,n}}\bigg)^{\frac{n}{2}} \nonumber \\
&\leq n \left( \frac{e}{\nu} \right)^{\nu n}  \left( 1-\mu_r+\mu_re^{- \mu_r} + e^{- \gamma K_{1,n}} \right)^{\frac{n}{2}}
\nonumber \\
&= n \left( \left( \frac{e}{\nu} \right)^{2\nu} \left( 1-\mu_r+\mu_re^{- \mu_r} + e^{-\gamma K_{1,n} }\right)\right)^{\frac{n}{2}}
\label{eq:conn_15}
\end{align}

We have $1-\mu_r+\mu_re^{- \mu_r}<1$ from $\mu_r>0$ and 
$e^{- \gamma K_{1,n}}=o(1)$ from
 (\ref{eq:conn_K1_new_trick_2}). Also,  it holds that $\lim_{\nu \to 0} \left( \frac{e}{\nu}\right)^{2\nu}=1$.
Thus, if we pick $\nu$ small enough to ensure that
\begin{equation}
\left( \frac{e}{\nu} \right)^{2\nu} \left( 1-\mu_r+\mu_re^{- \mu_r} \right) < 1,
\label{eq:conn_16}
\end{equation}
then for any $0< \epsilon < 1- \left( {e}/{\nu} \right)^{2\nu} \left( 1-\mu_r+\mu_re^{- \mu_r} \right)$ there exists a finite integer $n^{\star}(\epsilon)$  such that
\[
 \left( \frac{e}{\nu} \right)^{2\nu} \left( 1-\mu_r+\mu_re^{- \mu_r} + e^{-\gamma K_{1,n} }\right) \leq 1- \epsilon, \quad  \forall n \geq n^{\star}(\epsilon).
\]
Reporting this into (\ref{eq:conn_15}), we get
\begin{equation} \nonumber
\lim_{n \to \infty} \sum_{\ell=\left \lfloor{\frac{\mu_r n}{2 \beta c \log n}}\right \rfloor+1}^{\left \lfloor{\nu n}\right \rfloor} \binom{n}{\ell} \mathbb{P}[A_{n,\ell}(\pmb{\mu},\pmb{\Theta}_n) \cap E_n(\pmb{\mu},\pmb{\theta}_n,\pmb{X}_n)^c]=0
\end{equation}
since $\lim_{n \to \infty} n (1-\epsilon)^{n/2} = 0$.

\subsubsection{The case where $\left \lfloor{\nu n}\right \rfloor+1 \leq \ell \leq \lfloor{\frac{n}{2}} \rfloor $}
\label{subsec:range_5}
In this range, we use 
 (\ref{eq:conn_bounds_2}), (\ref{eq:conn_key_bound}),  the first bound in (\ref{eq:conn_key_bound_aux1}), the last bound in (\ref{eq:conn_key_bound_aux2}), and (\ref{eq:bounding_n_minus_l}) to get
 \begin{align*}
&\sum_{\ell=\left \lfloor{\nu n}\right \rfloor+1}^{\left \lfloor{\frac{n}{2}}\right \rfloor} \binom{n}{\ell} \mathbb{P}[A_{n,\ell}(\pmb{\mu},\pmb{\Theta}_n) \cap E_n(\pmb{\mu},\pmb{\theta}_n,\pmb{X}_n)^c] \nonumber \\
&\leq \sum_{\ell=\left \lfloor{\nu n}\right \rfloor+1}^{\left \lfloor{\frac{n}{2}}\right \rfloor} \binom{n}{\ell} \left( e^{- \beta \ell \alpha_n p_{11}(n)} + e^{- \gamma K_{1,n}}\right)^{\frac{n}{2}}\nonumber \\
&\leq \left(\sum_{\ell=\left \lfloor{\nu n}\right \rfloor+1}^{\left \lfloor{\frac{n}{2}}\right \rfloor} \binom{n}{\ell} \right) \left( e^{- \beta \nu n \alpha_n p_{11}(n)} + e^{- \gamma K_{1,n}}\right)^{\frac{n}{2}} \nonumber
\\
& \leq  \left( 4 e^{- \beta \nu n \alpha_n p_{11}(n)} + 4 e^{- \gamma K_{1,n}}\right)^{\frac{n}{2}} 
\end{align*}
 
With $\beta,\nu,\gamma>0$ have
$e^{- \beta \nu n \alpha_n p_{11}(n)} = o(1)$ from (\ref{eq:conn_K1}) 
and $e^{- \gamma K_{1,n}}=o(1)$ from
 (\ref{eq:conn_K1_new_trick_2}).
The conclusion
\begin{equation} \nonumber
\lim_{n \to \infty} \sum_{\ell=\left \lfloor{\nu n}\right \rfloor+1}^{\left \lfloor{\frac{n}{2}}\right \rfloor} \binom{n}{\ell} \mathbb{P}[A_{n,\ell}(\pmb{\mu},\pmb{\Theta}_n) \cap E_n(\pmb{\mu},\pmb{\theta}_n,\pmb{X}_n)^c]=0
\end{equation}
immediately follows and
 the proof of  one-law is completed.
\myendpf

\end{document}